\documentclass[12pt]{article}
\usepackage{mathrsfs}
\usepackage{graphicx}
\usepackage{amsmath}
\usepackage{amssymb}
\usepackage{amsthm}

\usepackage{thmtools}
\usepackage{amsfonts}

\usepackage[utf8]{inputenc}
\usepackage{tikz-cd}
\usetikzlibrary{tikzmark,quotes}

\usepackage{bbding}
\usepackage{comment}
\usepackage{mathtools}
\usepackage{xcolor}

\usepackage{ltablex}
\usepackage{float}
\usepackage[nospace,noadjust]{cite}
\usepackage{calc}
\usepackage{geometry}

\usepackage{indentfirst}
\usepackage{longtable}

\allowdisplaybreaks

\numberwithin{equation}{section} %\usepackage{geometry}

\DeclareMathOperator{\Supp}{Supp}
\DeclareMathOperator{\Exch}{Exc_{hol}}
\DeclareMathOperator{\exc}{exc}

\newtheorem{definition}{Definition}[section]

\newtheorem{theorem}[definition]{Theorem}

\newtheorem{corollary}[definition]{Corollary}
\newtheorem{lemma}[definition]{Lemma}
\newtheorem{conjecture}[definition]{Conjecture}

\newtheorem{example}[definition]{Example}
\newtheorem*{remark}{Remark}
\newtheorem{numberedremark}[definition]{Remark}

\newcommand{\Qbar}{ \overline{\mathbb{Q}}}
\newcommand{\Kbar}{ \overline{k}}

\newcommand{\bP}{\mathbb{P}}

\newcommand{\NS}{\mathrm{NS}}
\newcommand{\Pic}{\mathrm{Pic}}

\newcommand{\Gm}{\mathbb{G}_m}
\newcommand{\Pn}{\mathbb{P}^n}

\DeclareMathOperator{\rk}{rk}
\DeclareMathOperator{\codim}{codim}
\newcommand{\vol}{\operatorname{vol}}

\newcommand{\Addresses}{{% additional braces for segregating \footnotesize
  \bigskip
  \footnotesize
  \noindent Keping Huang, \textsc{Institute for Advanced Studies in Mathematics, Harbin Institute of Technology, Harbin, China 150001}\\
    \noindent \textit{E-mail address}: \texttt{kphuang@hit.edu.cn}
\bigbreak
  \noindent Aaron Levin, \textsc{Department of Mathematics, Michigan State University, East Lansing, USA 48824}\\
    \noindent \textit{E-mail address}: \texttt{levina@msu.edu}

}
}

\usepackage{hyperref}
\hypersetup{
    colorlinks=true,
    linkcolor=blue,
    filecolor=blue,      
    urlcolor=cyan,
}

\title {Greatest Common Divisors and an Unlikely Intersection Type Conjecture}

\author{Keping Huang and Aaron Levin}

\date{\vspace{-2em} }

\begin{document}
\maketitle

\def\Z{{\bf Z}}
%{Split Case, over $ \mathbb{C}$}  \hfill Authors
%\vspace{1em}

%%PLEASE PUT YOUR ABSTRACT HERE
\begin{abstract}
We prove inequalities involving greatest common divisors of functions at integral points with respect to numerically parallel divisors, generalizing a result of Wang and Yasufuku (after work of Bugeaud-Corvaja-Zannier, Corvaja-Zannier, and the second author). In addition to proving results in a broader context than previously considered, we also study the exceptional set in this setting, for both the counting function and the proximity function. Based on our study of exceptional sets, we propose new conjectures connected to the topic of unlikely intersections. We prove special cases of the conjectures on abelian varieties and the algebraic torus, including a completely explicit result in the former case, and provide a link between the conjectures and exceptional sets for inequalities involving greatest common divisors.
%After applying a result of Vojta on integral points on subvarieties of semiabelian varieties, we use geometry and the theory of heights to reduce to the (known) case of $\Gm^n$. 
\end{abstract}
%%THE END OF ABSTRACT

\section{Introduction}

We study the following three closely related objects: height inequalities which may be interpreted as upper bounds for certain (generalized) greatest common divisors, exceptional sets in these height inequalities, and a purely geometric conjecture which may be viewed as belonging to the subject of ``unlikely intersections".

\subsection{An Unlikely Intersection Conjecture}

We begin by stating a version of our geometric conjecture, both since it may be of independent interest and because it requires little background.

\begin{conjecture}\label{conj: unlikely3}
Let $X$ be a nonsingular complex projective variety with nonnegative Kodaira dimension $\kappa(X)\geq 0$. Let $Y$ be a positive-dimensional subvariety of $X$ of codimension at least $2$. Let $\mathbf{W}$ be the union of the subvarieties $W$ of $X$ such that
\begin{enumerate}
\item $W$ is not a variety of general type.\label{mainconj1}
\item $Y\subset W$ and $Y$ is an ample Cartier divisor on $W$.\label{mainconj2}
\end{enumerate}
Then $\mathbf{W}$ is not Zariski dense in $X$.
\end{conjecture}

More generally, we formulate a version of the conjecture when the log Kodaira dimension of $X$ is nonnegative (so $X$ is not required to be compact). 

We view the condition that $Y\subset W$ is an ample divisor on $W$ as being an ``unlikely intersection" type condition. Indeed, when $X=A$ is an abelian variety, we make a direct connection to classical unlikely intersection problems. Let $W$ be a subvariety of $A$ satisfying \eqref{mainconj1} and \eqref{mainconj2} (and $\dim Y>0$). Then since $W$ is not of general type, by a result of Ueno, $W$ can be written as a union of cosets of the form $a+B$, where $B$ is a positive-dimensional abelian subvariety of $A$ and $a\in W(\mathbb{C})$. If $\dim B\geq 2$ and $a+B\subset W$, then since $Y$ is ample in $W$ and $Y$ has codimension at least $2$ in $A$, 
\begin{align*}
\dim Y\cap (a+B)\geq \dim B-1>\dim Y+\dim(a+B)-\dim A, 
\end{align*}
and for each such coset we obtain a classical unlikely intersection. If $\dim B=1$ and $a+B\subset W$, then one can show (Theorem \ref{comparisontheorem}) that for some $\delta>0$ (independent of $W$ and $B$), 
\begin{align}
\label{unlikelyineq}
\deg(Y\cap (a+B))\geq \delta \deg B, 
\end{align}
where we define $\deg B$ by fixing some ample divisor on $A$ (or fixing a projective embedding of $A$). Again, an intersection satisfying \eqref{unlikelyineq} may be viewed as unlikely; see the discussion in work of the second author and Corvaja and Zannier \cite{CLZ22}, where analogous intersections in the case of the algebraic torus were studied and applied to problems involving lacunary polynomials and to other unlikely intersection problems. In this context, Conjecture \ref{conj: unlikely3} can be compared with the main Structure Theorem in the foundational work of Bombieri, Masser, and Zannier \cite[Theorem 1.4]{BMZ} asserting (notation as in \cite{BMZ}) that $H\mathcal{Z}_H, H\in \Phi$, is not Zariski dense in $\mathbb{G}_m^n$. From this viewpoint, our conjecture replaces abelian (or semiabelian) varieties by the much broader class of varieties of nonnegative log Kodaira dimension, and replaces translates of algebraic subgroups by subvarieties which are not of log general type; to compensate for this weakening, the conjecture requires a stronger unlikely intersection condition ($Y\subset W$ ample on $W$).

We now briefly discuss the relation between Conjecture \ref{conj: unlikely3} and inequalities for greatest common divisors. That there is some relation between unlikely intersections and greatest common divisor inequalities is not new; as noted by Zannier \cite[Chapter II]{Zan12}, Theorem \ref{tBCZ} below and its generalizations may be viewed as arithmetic analogues of unlikely intersection problems. Moreover, the aforementioned paper \cite{CLZ22} directly connected function field versions of GCD inequalities to unlikely intersection problems. Our motivation for Conjecture \ref{conj: unlikely3} originated in analyzing the form of the exceptional set in the conclusion of Theorem~\ref{thm: exceptional}, and our attempts to understand the converse question of when a closed subset must appear in such an exceptional set. This line of thinking culminated in Theorem \ref{thm: unlikelyexceptional}, which proves Conjecture \ref{conj: unlikely2}, conditional on well-known conjectures of Lang and Vojta, and provides a formal connection with exceptional sets in GCD type inequalities. This connection also provides additional insight into the conjectures; for instance, the necessary and subtle intersection condition on $Y$ and $D$ in Conjecture \ref{conj: unlikely2} may have been overlooked otherwise (see Example \ref{intersectionexample}).

Overall, we provide two primary pieces of evidence for Conjecture \ref{conj: unlikely3}, and the more general Conjecture \ref{conj: unlikely2}:
\begin{enumerate}
\item We prove that Conjecture \ref{conj: unlikely2} follows from conjectures of Lang and Griffiths/Vojta (Theorem \ref{thm: unlikelyexceptional}). 
\item We prove Conjecture \ref{conj: unlikely2} when $X$ is an abelian variety or a smooth toric compactification of $\mathbb{G}_m^n$ (Theorem \ref{abelian2} and Corollary \ref{toriccorollary}).
%\item We prove Conjecture \ref{conj: unlikely1} when $X$ is an abelian variety and $\dim Y < \frac{1}{3}\dim X$ (Theorem \ref{abelian1}).
\end{enumerate} 
In fact, when $X$ is an abelian variety we prove a completely explicit result (Theorem~\ref{abelian2}), which may be of independent interest.

Aside from intrinsic arithmetic interest, we believe the close connection between these conjectures and GCD inequalities provides additional geometric motivation to study GCD inequalities in the broadest possible context, with a focus on the exceptional set appearing in the inequalities. Indeed, the form of Theorem~\ref{thm: exceptional} provided the original motivation for formulating the conjectures given here.

\subsection{GCD Inequalities}

In 2003, Bugeaud, Corvaja, and Zannier \cite{BCZ} began a new line of results with their proof of the following inequality involving greatest common divisors:

\begin{theorem}[Bugeaud, Corvaja, Zannier \cite{BCZ}]
\label{tBCZ}
Let $a,b\in \mathbb{Z}$ be multiplicatively independent integers.  Then for every $\varepsilon>0$, 
\begin{align}
\label{eBCZ}
\log \gcd (a^n-1,b^n-1)\leq \varepsilon n
\end{align}
for all but finitely many positive integers $n$.
\end{theorem}

The inequality \eqref{eBCZ} was subsequently generalized by Corvaja and Zannier \cite{CZ05}, allowing $a^n$ and $b^n$ to be replaced by elements $u$ and $v$, respectively, of a fixed finitely generated subgroup of $\Qbar^*$, and replacing  $u-1$ and $v-1$ by more general pairs of polynomials in $u$ and $v$. The second author \cite{Lev19} further generalized this inequality to multivariate polynomials, resulting in the following statement: 

\begin{theorem}
\label{Lgcd}
Let $n$ be a positive integer, $\Gamma\subset \mathbb{G}_m^n(\Qbar)$ a finitely generated group, and $f(x_1,\ldots, x_n), g(x_1,\ldots, x_n)\in \Qbar[x_1,\ldots, x_n]$ nonconstant coprime polynomials such that not both of them vanish at $(0,0,\ldots,0)$.  Let $h(\alpha)$ denote the (absolute logarithmic) height of an algebraic number $\alpha$.  For all $\varepsilon>0$, there exists a finite union $Z$ of translates of proper algebraic subgroups of $\mathbb{G}_m^n$ such that
\begin{align}
\label{eLgcd}
\log \gcd (f(u_1,\ldots, u_n),g(u_1,\ldots, u_n))<\varepsilon \max\{h(u_1),\ldots, h(u_n)\}
\end{align}
for all $(u_1,\ldots, u_n)\in \Gamma\setminus Z$.
\end{theorem}

Silverman \cite{Sil05} interpreted (\ref{eBCZ}) as a special case of Vojta's Conjecture, and derived a more general conjectural inequality
containing (\ref{eLgcd}) as a special case.

If we identify $ \mathbb{G}_m^n$ with $\mathbb{P}^n\setminus \bigcup_{i=0}^{n}H_i$, where $H_0,\ldots, H_n$ are the $n+1$ coordinate hyperplanes, and let $Y$ be the closed subscheme of $\mathbb{P}^n$ defined by the ideal generated by the (homogenizations) of $f$ and $g$, then the height $h_Y$ associated to $Y$ is closely related to the left-hand side of \eqref{eLgcd}.  In this language, the following counterpart to Theorem \ref{Lgcd} was proven in \cite{Lev19}.

\begin{theorem}
\label{thpn}
Let $Y$ be a closed subscheme of $\mathbb{P}^n$, defined over a number field, of codimension at least $2$.  Suppose that $Y$ does not contain any of the points $P_0=[1:0:\cdots :0], \ldots, P_n=[0:0:\cdots:0 :1]$. Let $\Gamma\subset\mathbb{G}_m^n(\Qbar)$ be a finitely generated subgroup.  Then for all $\varepsilon>0$, there exists a finite union $Z$ of translates of proper algebraic subgroups of $\mathbb{G}_m^n$ such that
\begin{align*}
h_Y(P)\leq \varepsilon h(P)
\end{align*}
for all $P\in \Gamma\setminus Z\subset \mathbb{P}^n(\Qbar)$.
\end{theorem}

\begin{sloppypar}
Following Wang and Yasufuku \cite{WY19}, we study extensions of Theorem \ref{thpn} where $\mathbb{P}^n$ is replaced by an arbitrary projective variety $X$, the coordinate hyperplanes $H_0,\dots, H_{n}$ are replaced by numerically equivalent (or more generally numerically parallel) divisors $D_1,\ldots, D_{n+1}$ on $X$, and the finitely generated subgroup $\Gamma$,
which can be viewed as a subset of a set of $(\sum_{i=0}^n H_i, S)$-integral points (for appropriate $S$), 
is replaced by a set $R$ of $(\sum_{i=1}^{n+1}D_i,S)$-integral points on $X$. More precisely, and to highlight the improvements and assumptions in our results, we work in the following {\bf general setup}:
\end{sloppypar}

\begin{center}
\begin{tabular}{p{2.5cm}p{11cm}}
$X$ &  a projective variety of dimension $n$ defined over a number field $k$.  \\
$\mathrm{Div}(X)$ & the group of Cartier divisors on $X$.\\
$D_1,\dots,D_{n+1}$   & effective Cartier divisors on $X$, defined over $k$, that are $\mathbb{Z}$-linearly independent in $\mathrm{Div}(X)$, and such that there exist positive integers $d_i$ such that $d_jD_i \equiv d_i D_j$ for all $1 \le i,j \le n + 1$. \\
$D$ & the sum $\sum_{i=1}^{n+1} D_i$.\\
$Y$ &   a closed subscheme of $X$, defined over $k$, of codimension at least $2$. \\
$M_k$ & the set of (nontrivial) places of $k$. \\
$S\subset M_k$ & a finite set of places of $k$ containing the archimedean places. \\
$\mathcal{O}_{k,S}$ & the ring of $S$-integers of $k$. \\
$R\subset X(k)$ & a set of  $(D, S)$-integral points on $X$. \\
\end{tabular}
\end{center}

\begin{comment}
\begin{itemize}
\item  $X$ is a projective variety of dimension $n$ defined over a number field $K$.
\item  $D_1,\dots,D_{n+1}$ are effective Cartier divisors of $X$, defined over $K$, that are linearly independent in $\mathrm{Div}(X)$, and there exist positive integers $d_i$ such that $d_jD_i \equiv d_i D_j$ for all $1 \le i,j \le n + 1$.
\item $Y$ is a closed subscheme of $X$, defined over $K$, of codimension at least $2$.
\item $S$ is a finite set of places of $K$ containing the archimedean places and $R\subset X(K)$ is a set of  $(D, S)$-integral points on $X$, where $D = \sum_{i=1}^{n+1} D_i$.
\end{itemize}
\end{comment}

Based on the powerful tools developed in \cite{RW20}, Wang and Yasufuku \cite{WY19} proved the following generalization of Theorem \ref{Lgcd}.

%theorem, which can be interpreted as an upper bound for the greatest common divisor for certain types of integral points. For example, it includes a main result of \cite{CZ05} as a special case. 

\begin{theorem}[Wang, Yasufuku \cite{WY19}]\label{thm: WY}
In addition to the general setup, suppose that $X$ is a Cohen–Macaulay variety, the divisors $D_1,\ldots, D_{n+1}$ are ample and in general position, and the closed subscheme $Y$ does not contain any point of the set
$\bigcup_{i=1}^{n+1}\bigcap_{j\neq i}\mathrm{Supp}D_j.$ 
Let $\varepsilon > 0$ and let $A$ be an ample divisor on $X$. Then there exists a proper Zariski closed subset $Z$ of $X$ such that
$$ h_Y (P) \le  \varepsilon h_A(P)$$
for all points $P \in R \setminus Z$. 
\end{theorem}

%Theorem \ref{thm: WY} is related with Theorem \ref{thpn} in the following sense. 

The main result of Section \ref{sec: pf} is the following theorem. 

\begin{theorem}\label{thm: height}
In addition to the general setup, suppose that the closed subscheme $Y$ does not contain any point of the set
$\bigcup_{i=1}^{n+1}\bigcap_{j\neq i}\mathrm{Supp}D_j.$ 
Let $\varepsilon > 0$ and let $A$ be a big divisor on $X$. Then there exists a proper Zariski closed subset $Z$ of $X$ such that
\begin{equation} \label{eq: main}
h_{Y}(P) < \varepsilon  h_A(P)
\end{equation}
for all $P \in R \setminus Z$. 
\end{theorem}

Note that the conditions on the variety $X$ and the divisors $D_i$ are weaker than those in Theorem \ref{thm: WY}. 
In particular, we do not require the divisors $D_i$ to be in general position, 
nor do we require the divisors $D_i$ to be ample. 
In fact, the issue of non-general position is important, since even if we are only interested in divisors in general position, when we consider inductive arguments to prove finiteness statements, one naturally encounters restrictions of the divisors to the exceptional set, and the resulting divisors may no longer be in general position.

The geometric condition on $Y$ in Theorems \ref{thpn}--\ref{thm: height} (or the nonvanishing condition in Theorem \ref{Lgcd}) cannot be dropped in general. However, going back to the work of Corvaja and Zannier \cite[Prop.~4]{CZ05}, it has been observed that this geometric condition may be dropped if one avoids a finite number of local contributions to $h_Y$ (or the logarithmic gcd).

%\begin{remark}
%In fact, inequality (\ref{eq: main}) is vacuous in the first case above, 
%so it suffices to conclude the second case, but for clarity, we divide into two cases. 
%\end{remark}

More precisely, outside of the support of $Y$, we can decompose the height $h_Y$ into a sum of local heights
\begin{align*}
h_Y(P)=\sum_{v\in M_k}\lambda_{Y,v}(P),
\end{align*}
and write $h_Y=m_{Y,S}+N_{Y,S}$, where $m_{Y,S}=\sum_{v\in S}\lambda_{Y,v}$ and $N_{Y,S}=\sum_{v\in M_k\setminus S}\lambda_{Y,v}$ are known as the proximity function and counting function, respectively, associated to $Y$ and the finite set of places $S\subset M_k$. For the counting function, the following result was proved in \cite{Lev19}.

\begin{theorem}
\label{Lgcd2}
Let $Y$ be a closed subscheme of $\mathbb{P}^n$ of codimension at least $2$, defined over a number field $k$.  Let $S$ be a finite set of places of $k$ containing the archimedean places.  For all $\varepsilon>0$, there exists a finite union $Z$ of translates of proper algebraic subgroups of $\mathbb{G}_m^n$ such that
\begin{align*}
N_{Y,S}(P)\leq \varepsilon h(P)
\end{align*}
for all $P\in \mathbb{G}_m^n(\mathcal{O}_{k,S})\setminus Z\subset \mathbb{P}^n(k)$.
\end{theorem}

Theorem \ref{Lgcd2} was generalized by Wang and Yasufuku under the assumption that $Y$ consists of a finite number of points:
\begin{theorem}[\cite{WY19}]\label{thm: WY2}
In addition to the general setup, suppose that $X$ is a Cohen–Macaulay variety, the divisors $D_1,\ldots, D_{n+1}$ are ample and in general position, and $\dim Y=0$.  Let $\varepsilon > 0$ and let $A$ be an ample divisor on $X$. Then there exists a proper Zariski closed subset $Z$ of $X$ such that
$$ N_{Y,S} (P) \le  \varepsilon h_A(P)$$
for all points $P \in R \setminus Z$. 
\end{theorem}

We give a complete generalization of Theorem \ref{Lgcd2} (and Theorem \ref{thm: WY2}) under the hypotheses of the general setup.

\begin{theorem}\label{thm: count}
Assume the general setup. Let $\varepsilon > 0$ and let $A$ be a big divisor on $X$. Then there exists a proper Zariski closed subset $Z$ of $X$ such that
\begin{equation} \label{eq: counting}
N_{Y,S}(P) < \varepsilon h_A(P)
\end{equation}
for all $P \in R \setminus Z$. 
\end{theorem}

\begin{comment}

The proof is divided into two parts, one for the counting function and the other for the proximity function. The following result on the counting function is of its own interest.

\begin{remark}
If in addition $A$ is ample, and there are $n+2$  effective Cartier divisors $D_i,~i=1,\dots, n+2$ such that $D_i\equiv d_i A,~d_i\in \mathbb{Z}_{\ge 0}$ and $D_i$ and $D_j$ don't share a common component $(i\neq  j)$, then we can apply the main theorem 
of \cite{RW20} and conclude the degeneracy of $(D,S)$-integral points. 
Therefore the number $n+1$ in Theorem \ref{thm: count} is optimal in this case to guarantee that we are in the second case. 
\end{remark}

Although it only proves a bound for the counting function, 
Theorem \ref{thm: count} weakens the condition that $D_i$ must be in general position, 
and eliminates the condition about the set $Y$ in Theorem \ref{thm: WY}.
We only require that $Y$ is of codimension at least $2$. 
In addition, we don't require that $A$ is ample. 
\end{comment}

In the case when the divisors $D_1,\dots,D_{n+1}$ are ample, 
we can say more about the exceptional $Z$ in Theorem \ref{thm: count}. We show that the exceptional set $Z$ may be chosen so that $Y$ intersects every irreducible component of $Z$ in a codimension one set.
\begin{theorem}\label{thm: exceptional}
In addition to the general setup, assume further that $\bigcap_{i=1}^{n+1}D_i = \emptyset$ and $D_1,\ldots, D_{n+1}$ are ample. Let $\varepsilon > 0$ and let $A$ be an ample divisor on $X$. Then there exists a proper Zariski closed subset $Z$ of $X$ such that for each irreducible component $Z'$ of $Z$, the subscheme $Z' \cap Y $ has codimension at most $1$ in $Z'$, and
\begin{equation} 
N_{Y,S}(P) < \varepsilon h_A(P)
\end{equation}
for all $P \in R \setminus Z$. 
\end{theorem}

The conclusion of Theorem \ref{thm: exceptional} naturally raises the problem of further studying the components of the exceptional set. Even if $Y\cap Z'$ is a divisor on $Z'$, the counting function (or the height $h_Y$) may still be small because the divisor is too ``small" (e.g., in an extreme case it could be an exceptional divisor for some morphism). Moreover, the set $Z'$ is not needed in the exceptional set if it lacks integral points. Taking these two further observations into account motivated Conjecture \ref{conj: unlikely3} and Conjecture \ref{conj: unlikely2}, and Theorem \ref{thm: unlikelyexceptional} formally connecting the latter conjecture to exceptional sets in GCD inequalities.

As a complement to Theorem \ref{thm: exceptional}, we similarly study the exceptional set for analogous inequalities where the counting function $N_{Y,S}$ is replaced by the proximity function $m_{Y,S}$ (Theorem \ref{thm: proxExc}).

%With Theorem \ref{thm: count}, we prove the following theorem, which is a more direct generalization of Theorem \ref{thm: WY}. 

%Theorem \ref{thm: height} generalize Theorem \ref{thm: WY} in two different perspectives. 

%Theorem \ref{thm: height} also doesn't require the assumption of general position as in Theorem \ref{thm: WY}. 

Our proof strategy for the results of this section differs from \cite{WY19}, where they mainly applied an inequality of Ru and Vojta \cite{RW20}. Since we don't assume the divisors are in general position, 
we cannot apply the Ru-Vojta result. 
Fortunately we are able to construct a morphism to $\mathbb{G}_m^n$ and reduce to the inequalities there, proved by the second author \cite{Lev19}. 
The proofs proceed in roughly two steps. First, we use a result of Vojta on integral points on subvarieties of semiabelian varieties to reduce to considering the case where the divisors $D_i$ are, up to some positive integer multiple, linearly equivalent (instead of just numerically equivalent). Second, we use the linear equivalence to construct a suitable morphism to projective space, where we apply functoriality of heights and results of the second author on $\mathbb{G}_m^n$ to derive the desired inequalities. In the second step, in the case of the proximity function, we prove the requisite result on $\mathbb{G}_m^n$ in Theorem \ref{thm: proximity}, and provide a more refined analysis in Section \ref{sec: toric}.

%The proof idea of the above theorems is to map $X$ to a compactification $\mathbb{P}^n$ of the algebraic torus $\Gm^n$ and use Theorems \ref{Lgcd}, \ref{Lgcd2} and other theorems by the second author \cite{Lev19} to derive the requested upper bounds. 
%We consider the %Albanese variety $\mathrm{Alb}(X)$ of $X$ and 
%irregularity $q=h^1(X,\mathcal{O}_X)$. 
%If $q > 0 $, we have the degeneracy of integral points by a main theorem of \cite{Voj96};  
%if $q=0$, we can apply Theorem \ref{Lgcd2} to prove Theorem \ref{thm: count}. 
%Combining Theorem \ref{thm: proximity} for proximity function and the assumption on the subscheme $Y$, 
%we are able to prove Theorem \ref{thm: height}. 

\subsection{Analogues in Value Distribution Theory}

Finally, we discuss analogues of our main results in Nevanlinna theory.  The proofs of our arithmetic results rely on gcd inequalities of the second author in \cite{Lev19}. 
Analogous inequalities in Nevanlinna theory, proved more generally for semiabelian varieties, are due to Noguchi, Winkelmann, and Yamanoi: 

\begin{theorem}[Noguchi, Winkelmann, Yamanoi \cite{NWY02}, Theorem 5.1]\label{thm: complex}
Let $f: \mathbb{C} \rightarrow A$ be a holomorphic map to a semiabelian variety $A$ 
with Zariski-dense image. Let $Y$ be a closed subscheme of $A$ with $\codim Y \ge 2$ and 
let $\varepsilon > 0$. 
\begin{enumerate}
    \item Then $$N_f(Y,r) \le_{\mathrm{exc}} \varepsilon T_f(r).$$
    \item There exists a compactification $\bar{A}$ of $A$, independent of $\varepsilon$, such that for the Zariski closure $\overline{Y}$ of $Y$ in $\overline{A}$, 
    $$T_{\overline{Y}, f}(r) \le_{\mathrm{exc}} \varepsilon T_f(r).$$
\end{enumerate}
\end{theorem}

Here $N_{f}(Y,r) $ is a counting function associated to $f$ and $Y$, $T_{\overline{Y}, f}(r)$ is a Nevanlinna characteristic function associated to $f$ and $\overline{Y}$, and $T_f(r)$ is any characteristic function associated to an appropriate ample line bundle. The notation $\le_{\mathrm{exc}} $ means that the estimate holds for all $r$ outside of a set of finite Lebesgue measure, possibly depending on $\varepsilon$. See also \cite{NW14}, Section 6.5 and \cite{LW20} for more discussion.  

We also need the analogue for holomorphic curves of Vojta's results on integral points on subvarieties of semiabelian varieties due to Noguchi \cite{Nog81} (as formulated by Noguchi-Winkelmann \cite[Th.~4.9.7]{NW14}).

Using these results, one can derive the following Nevanlinna-theoretic inequalities, analogous to Theorems \ref{thm: height} and \ref{thm: count}, and generalizing a result of Wang-Yasufuku \cite[Th.~1.8]{WY19}. We omit the details here. 

\begin{theorem}
Assume the general setup, with the number field $k$ replaced by $k=\mathbb{C}$. 
Let $\varepsilon > 0$ and let $A$ be a big divisor on $X$. 
Let $f: \mathbb{C} \rightarrow X \setminus \bigcup_{i=1}^{n+1} D_i$ be a holomorphic map with Zariski-dense image. Then we have 
$$ N_{f}(Y,r) \le_{\mathrm{exc}} \varepsilon T_{A,f}(r). $$
Suppose furthermore that $Y$ does not contain any point of the set $\bigcup_{i=1}^{n+1}\bigcap_{j\neq i} D_j$. Then we have 
$$ T_{Y,f}(r) \le_{\mathrm{exc}} \varepsilon T_{A,f}(r). $$
\end{theorem}

\subsection{Outline of the Paper}

The background and notation for our results will be given in Section \ref{sec: notation}. We prove Theorems \ref{thm: height} and \ref{thm: count} in Section \ref{sec: pf}. In Section \ref{unlikely} we study Conjecture \ref{conj: unlikely3} and Conjecture \ref{conj: unlikely2}, proving some special cases and showing Conjecture \ref{conj: unlikely2} follows from standard conjectures (Theorem \ref{thm: unlikelyexceptional}). Finally, Section \ref{sec: toric} is devoted to inequalities on toric varieties and their applications to the exceptional set for the proximity function under our general setup.

\section{Background and Notation}\label{sec: notation}

In \cite{Sil87}, Silverman generalized the Weil height machine for Cartier divisors to height functions on projective varieties with respect to closed subschemes. 
More precisely, let $X$ be a projective variety over a number field $k$, and let $Z(X)$ denote the set of closed subschemes of $X$. 
Note that the closed subschemes $Y \in Z(X)$ are in one-to-one correspondence with quasi-coherent ideal sheaves 
$\mathcal{I}_Y \subseteq \mathcal{O}_X$, and we identify a closed subscheme $Y$ with its ideal
sheaf $\mathcal{I}_{Y}$. Silverman assigned to each $Y \in Z(X)$ and each place $v\in M_k$ a local height function $\lambda_{Y,v}$, 
and to each $Y \in Z(X)$ a global height function $
h_Y=\sum_{v\in M_k}\lambda_{Y,v}$ that generalizes the Weil height machine for Cartier divisors.

\begin{theorem}[\cite{Sil87}]\label{thm: functoriality}
Let $X$ be a projective variety over a number field $k$. There are maps
\begin{equation*}
    \begin{aligned}
\lambda:  Z(X) \times M_k &\rightarrow \{\text{functions }X(\overline{k}) \rightarrow [0,+\infty] \},\\
h: Z(X) &\rightarrow \{\text{functions }X(\overline{k}) \rightarrow [0,+\infty] \},
    \end{aligned}
\end{equation*}
satisfying the following properties:
\begin{enumerate}
    \item   If $W,Y \in  Z(X)$ satisfy $W \subseteq Y$, then $h_W \le h_Y+O(1)$ and $\lambda_{W,v} \le \lambda_{Y,v}+O(1)$ for all $v\in M_k$. 
    \item  If $ W,Y \in Z(X)$ satisfy $\mathrm{Supp}(W) \subseteq \mathrm{Supp}(Y)$,
    then there exists a constant $C$  
such that $h_W \le C \cdot h_Y+O(1)$ and $\lambda_{W,v} \le C \cdot \lambda_{Y,v}+O(1)$ for all $v\in M_k$. 
		\item For all $W,Y\in Z(X)$, we have $h_{W+Y}=h_W+h_Y+O(1)$ and $\lambda_{W+Y,v}=\lambda_{W,v}+\lambda_{Y,v}+O(1)$ for all $v\in M_k$.
    \item Let $\phi: X' \rightarrow X$ be a morphism of projective varieties over $k$, and let $Y \in Z(X)$. Then
		\begin{align*}
		h_{X',\phi^*Y} &= h_{X, Y}\circ \phi+O(1),\\
		\lambda_{X', \phi^*Y,v} &= \lambda_{X, Y,v}\circ \phi+O(1),
		\end{align*}
		for all $v\in M_k$. 
\end{enumerate}
\end{theorem}

Here, $W\subset Y$, $W+Y$, and $\phi^*Y$ are all defined in terms of the associated ideal sheaves (see \cite{Sil87}).

\begin{remark}
Note that $\phi^*Y$ is the closed subscheme corresponding to the inverse image ideal sheaf 
$$\phi^{-1}\mathcal{I}_Y \cdot \mathcal{O}_{X'}$$
instead of the pullback of the ideal sheaf. See Caution 7.12.2 of \cite{Har77} for more details. 
\end{remark}

We shall reduce the estimates of counting functions to those over a compactification $\Pn$ of $\Gm^n$. 
For this purpose we need the following elementary fact (see, e.g., \cite[p.~112]{GW}). 
%\begin{lemma}\label{lem: pullback}
%Let $\phi: X  \rightarrow V$ be a morphism of projective varieties and let $Y \in Z(X)$. Then we have 
%$\phi^{-1}(\mathcal{I}_{\phi(Y)}) \cdot \mathcal{O}_X \subseteq \mathcal{I}_Y$. 
%%$Y\subseteq Z(\phi^{-1}(\mathcal{I}_{\phi(Y)}))$ as subschemes of $X$. 
%\end{lemma}
%
%\begin{proof}
%%It suffices to prove that we have a sheaf inclusion 
%We will prove it at the level of stalks. 
%Suppose $P\in X$ is any point and let $Q = \phi(P)$.
%Set $B = \mathcal{O}_{P,X}$, $A = \mathcal{O}_{Q,Y}$ and $I = \mathcal{I}_{P}$. 
%Now the inclusion is clear as we have $\phi_P^{-1}(I) \cdot B \subseteq I\cdot B$. 
%\end{proof}

\begin{lemma}\label{lem: pullback2}
Let $\phi: X  \rightarrow V$ be a morphism of projective varieties and let $Y \in Z(V)$. Then $\Supp \phi^*Y=\phi^{-1}(\Supp Y)$. 
\end{lemma}

We use the following definition of integrality (see \cite{Voj87} for more details). 
\begin{definition}\label{def: integral}
Let $D$ be an
effective Cartier divisor on $X$, let $R$ be a subset of 
$X(k) \setminus \Supp D$, and let $S$ be a finite set of places of $k$ containing all the archimedean places. 
We say that $R$ is a set of \emph{$(D, S)$-integral points} on $X$ if
there exists a global Weil function $\lambda_{D,v}$ and constants $c_v$, $v \in M_k$, such that 
\begin{itemize}
    \item  $c_v = 0$ for all but finitely many $v$, and 
    \item for all $v \in M_k\setminus S$, we have $\lambda_{D,v}(P) \le  c_v$ for all $P \in R$. 
\end{itemize}
\end{definition}

According to this definition, a finitely generated subgroup $\Gamma \subset \mathbb{G}_m^n(\mathcal{O}_{k,S})$ can be viewed as a subset of a set of $(\sum_{i=0}^n H_i, S)$-integral points, 
since we can choose $c_v = 0$ for all $v \in M_k \setminus S$. 
Therefore Theorems \ref{thm: WY} and \ref{thm: height} are related with 
Theorem \ref{thpn} (and hence with Theorems \ref{Lgcd} and \ref{tBCZ}).

We denote numerical and linear equivalence of divisors by $\equiv$ and $\sim$, respectively.

\section{Proofs of Theorems \ref{thm: height}, \ref{thm: count}, and \ref{thm: exceptional}}\label{sec: pf}

\begin{comment}
If $g>0$, then we apply Vojta's Corollary 0.3 in \cite{Voj96}, integral points on the complement of $D_1,...,D_{n+1}$ will not be Zariski dense (in this case we can take $\rho=1$ from the proof of that Corollary, by the numerical equivalence condition). 
\end{comment}

We need the following lemma, whose essential ingredient is a result of Vojta \cite{Voj96} on integral points on subvarieties of semiabelian varieties.

\begin{lemma}
\label{lem: reduction}
Assume the general setup. Suppose that the set of $(D,S)$-integral points $R$ is Zariski dense in $X$. Then $q=h^1(X,\mathcal{O}_X)=0$ and for some positive integer $m$,  $md_jD_i \sim md_i D_j$ (linear equivalence) for all $i$ and $j$.  If $\mathrm{div}(f_i)=m(d_{n+1}D_i-d_iD_{n+1})$, $i=1,\ldots, n$, then the morphism
\begin{align*}
\psi&: X\setminus D \rightarrow \Gm^n,\\
\psi(P) &= \left({f_1}(P), \dots, {f_n}(P)\right)
\end{align*}
is dominant and there exists a finite set of places $S'\supset S$ such that $\psi(R)\subset \mathbb{G}_m^n(\mathcal{O}_{k,S'})$.
\end{lemma}

\begin{proof}
Suppose first that $q>0$. Then $D$ contains the $n+1$ linearly independent divisors $D_i$ and $n+1\geq \dim X-q+\rk \{D_i\}_{i=1}^{n+1}+1=n+2-q$, where $\rk \{D_i\}_{i=1}^{n+1}$ denotes the rank of the subgroup generated by $D_1,\ldots, D_{n+1}$ in the N\'eron-Severi group $\NS(X)$. Then $R$ is not Zariski dense by a result of Vojta \cite[Corollary 0.3]{Voj96} (noting that the Picard number $\rho$ in Vojta's statement may clearly be replaced by $\rk \{D_i\}_{i=1}^{n+1}$ in the proof).

Suppose now that $q=0$. By \cite[Corollary 9.5.14]{Klei05}, this implies $\Pic^0(X)$ is trivial. From the exact sequence (of groups)
\begin{align*}
0 \rightarrow \Pic^0(X) \rightarrow \mathrm{Pic}(X) \rightarrow \NS(X) \rightarrow 0,
\end{align*}
we have $\Pic(X) \cong \NS(X)$. Since by assumption $d_jD_i \equiv d_i D_j$ for all $i$ and $j$, and numerical triviality implies torsion in $\NS(X)$, it follows that for some positive integer $m$, $md_jD_i \sim md_i D_j$ for all $i$ and $j$. Let $f_1,\dots, f_{n} \in k(X)$ be rational functions such that $\mathrm{div}(f_i) = m(d_{n+1}D_i - d_iD_{n+1})$, $i=1,\ldots, n$. Define a morphism
$\psi: X\setminus D \rightarrow \Gm^n$ by $\psi(P) = \left({f_1}(P), \dots, {f_n}(P)\right)$.
Since $R\subset X(k)$ is a set of $(D,S)$-integral points and the functions $f_i$ have no zeros or poles outside $D$, by \cite[Lemma~1.4.6]{Voj87} there exists a constant $c\in k^*$ such that $c f_i(P), c/f_i(P) \in \mathcal{O}_{k,S} $ for all $P\in R$. Then for some finite set of places $S'\supset S$, 
we have 
\begin{equation}\label{eq: unit}
 f_i(P) \in \mathcal{O}_{k,S'}^\times
\end{equation}
for all $P\in R$, where $\mathcal{O}_{k,S'}^\times$ is the group of $S'$-units of $k$. It follows that $\psi(R)\subset \mathbb{G}_m^n(\mathcal{O}_{k,S'})$.

If $\psi$ is not dominant, then $\psi(R)$ is not Zariski dense in $\mathbb{G}_m^n$, and by Laurent's theorem, $\psi(R)$ is contained in a finite union $T$ of translates of proper algebraic subgroups of $\Gm^n$. Since $R$ is Zariski dense in $X$, it follows that $\psi(X)\subset T$.  But $X$ is irreducible, and so $\psi(X)$ is in fact contained in an irreducible component of $T$.  It follows that there is a nontrivial multiplicative relation
$f_1^{e_1}\dots f_n^{e_n}  = c$ with $c\in \bar{k}^\times$. 
Since $\mathrm{div}(f_i) = m(d_{n+1}D_i-d_iD_{n+1})$, this yields a nontrivial linear dependence relation for  
$D_1,\dots, D_{n+1}$, contradicting our {\bf general setup}. It follows that $\psi$ is dominant. 
\end{proof}

%Let $q=h^1(X,\mathcal{O}_X)$. 
%If $q>0$, we follow the proof of Corollary 0.3 in \cite{Voj96} to conclude the degeneracy of $(D,S)$-integral points. 
%In fact, replacing $D_i$ by an integral multiple 
%(which does not affect integrality), 
%we may assume that the divisors $D_i$ are all numerically equivalent. 
%Consider the divisors 
%$$F_i := D_i -  D_{n+1}\equiv 0,~i = 1,\dots, n. $$ 
%%$$F_i := d_{n+1} D_i - d_i D_{n+1}\equiv 0,~i = 1,\dots, n. $$ 
%By the work of Matsusaka \cite{Mat57}, 
%algebraic equivalence and numerical equivalence coincide for divisors up to torsion. So 
%some multiples $m F_i~(m\in \mathbb{Z}_{>0})$ are simultaneously algebraically equivalent to $0$. 
%The divisors $mF_i$ can be used to define a semiabelian variety $A$, 
%which $X \setminus D$ maps into 
%(see Section 4.5 of \cite{NW14}, \cite{NW02} or \cite{Iit76} for more details). 
%Then we have an exact sequence 
%$$ 1 \rightarrow  \Gm^{n} \rightarrow  A \rightarrow  \mathrm{Alb}(X) \rightarrow  1. $$
%%But by the exact sequence (2.3) in \cite{NW02}, we have
%It follows that
%$$\mathrm{dim} A = n + q > \mathrm{dim} X, $$
%and hence the image of $X \setminus D$ is a proper subvariety, to which we can apply Theorem 0.2 of \cite{Voj96}.
%This implies the degeneracy of any set of $(D,S)$-integral points on $X$. 

We now give the proof of Theorem \ref{thm: count}.

\begin{proof}[Proof of Theorem \ref{thm: count}]

First, by analyzing the irreducible components of $Y$ one by one, we may assume that $Y$ is itself irreducible.  By part 2 of Theorem \ref{thm: functoriality} we may assume further that $Y$ is reduced. Finally, we may assume that $R$ is Zariski dense in $X$.

Let $\psi:X\dashrightarrow \mathbb{P}^n$ be the rational map defined in Lemma \ref{lem: reduction}.  By blowing up the indeterminacy locus of $\psi$, we may find a projective variety $X'$ and morphisms $\pi:X'\to X$ and $\phi:X'\to \mathbb{P}^n$ such that $\pi$ is birational and $\psi\circ \pi=\phi$ on $\pi^{-1}(X\setminus D)$.

\begin{figure}[h!]
\centering
\begin{comment}
\begin{tikzcd}
   \subnode{u1}{X'} \arrow[d, ""'] 
  \arrow[rd, "\phi"]  &\\
 %\subnode{u4}{Y}  \subset 
 \subnode{u3}{X}  \arrow[r, "\psi"', dotted] & \mathbb{P}^n \supset W
\end{tikzcd}
\begin{tikzpicture}[overlay, remember picture]
  \draw [->] (u1) edge ["$\pi$", right] (u3.north -| u1);
\end{tikzpicture}
\end{comment}

% https://tikzcd.yichuanshen.de/#N4Igdg9gJgpgziAXAbVABwnAlgFyxMJZARgBoAGAXVJADcBDAGwFcYkQANAchAF9T0mXPkIoyxanSat2HPgJAZseAkQBMpCTQYs2iEAB0DAW3o4AFgCNLwAAq8AeoX6DlIogGZNknTP0B1eVdhVRRyb21pPRAATT5JGCgAc3giUAAzACcIYyRwkBwIJGIXECycvJpCpDVS8tzEMgKixA0QRiwwaKgIHBxEoLLshq9mmqr6LEZ2cwgIAGtB+qQAFiqWppxJ6f1ZhfjeIA
\begin{tikzcd}
         & X' \arrow[d, "\pi"] \arrow[rd, "\phi"] &              &                   \\
Y \arrow[r, hookrightarrow] & X \arrow[r, "\psi", dotted]   & \mathbb{P}^n \arrow[r, hookleftarrow] & W% \arrow[l,hookrightarrow
\end{tikzcd}
\caption{The Geometric Setup}
\end{figure}

By Definition \ref{def: integral}, for all $P\in R$ and $P'=\pi^{-1}(P)$,
\begin{align}
\label{NDeq}
N_{D,S}(P)=N_{\pi^*D,S}(P')+O(1)\leq O(1).
\end{align}

Suppose first that $Y\subseteq D$. In this case, by \eqref{NDeq}, for all $P\in R$,  
we have $$N_{Y,S}(P) \le N_{D,S}(P) < C,$$
where $C $ depends only on $R$. The result now follows from an appropriate version of Northcott's theorem \cite[Prop.~1.2.9(h)]{Voj87}.

Now assume that $Y \not \subseteq D$.  Therefore $\Supp(Y\cap D)$ is a proper Zariski closed subset of $Y$. Let $W$ be the closed subscheme with the reduced induced structure on the Zariski closure of $\phi(\pi^{-1}(Y\setminus D))=\psi(Y\setminus D)$. Then $\dim W=\dim(\psi(Y\setminus D))\leq \dim(Y\setminus D)=\dim Y$.   Since $\dim X=n$ and $\codim Y\geq 2$, it follows that $W$ is of codimension at least $2$ in $\mathbb{P}^n$.

By Theorem \ref{Lgcd2}, for any given $\varepsilon > 0 $, 
there exists a proper Zariski closed subset $Z_0\subset\mathbb{P}^n$ 
such that for all $\psi(P)  \in \psi(R)\setminus Z_0$, we have 
$$N_{W,S}(\psi(P)) \le \varepsilon \cdot h(\psi(P)) + O(1)$$
(strictly speaking, this only follows from Theorem \ref{Lgcd2} with $S$ replaced by the $S'$ of Lemma \ref{lem: reduction}; however, Theorem \ref{Lgcd2} easily extends to sets of $(D,S)$-integral points, where $D$ is the boundary of $\mathbb{G}_m^n\subset\mathbb{P}^n$, and $\psi(R)$ is a set of $(D,S)$-integral points from the proof of Lemma \ref{lem: reduction}).
Equivalently, if $P' = \pi^{-1}(P)$, which is well-defined as $\pi$ is an isomorphism over $X\setminus D$, then
$$N_{W,S}(\phi(P'))  \le \varepsilon \cdot h(\phi(P')) + O(1). $$
If $A'=\phi^*H$, where $H$ is a hyperplane of $\mathbb{P}^n$, then $A'$ is a big divisor and by Theorem \ref{thm: functoriality},  
 for all $P' \in \pi^{-1}(R)$ such that $\phi(P')\notin Z_0$, we have
\begin{align}
\label{NWeq}
N_{\phi^*W,S}(P') \le \varepsilon   h_{A'}(P') + O(1).
\end{align}
Now we note that
\begin{align*}
\pi^{-1}(Y)\subseteq\pi^{-1}(Y\setminus D)\cup \pi^{-1}(D)\subseteq \phi^{-1}(\psi(Y\setminus D))\cup \pi^{-1}(D)\subseteq \phi^{-1}(W)\cup \pi^{-1}(D),
\end{align*}
and it follows from Lemma \ref{lem: pullback2} that $\Supp \pi^*Y\subseteq \Supp(\phi^*W+\pi^*D)$. By Theorem \ref{thm: functoriality}, there exists a constant $C$ such that
\begin{align*}
N_{\pi^*Y,S}(P')\leq C(N_{\phi^*W,S}(P')+N_{\pi^*D,S}(P'))+O(1).
\end{align*}

Then by \eqref{NDeq} and \eqref{NWeq}, for all $P' \in \pi^{-1}(R)\setminus \phi^{-1}(Z_0)$, we find that
\begin{align*}
N_{\pi^*Y,S}(P')\leq \varepsilon h_{A'}(P') + O(1).
\end{align*}

Recalling that $\pi$ is birational and $\psi$ is dominant, combining this with Lemma 2.2 of \cite{Lev19}, there exists a constant $C' > 0$ and a proper Zariski closed subset $Z\supset \pi(\phi^{-1}(Z_0))$ of $X$ such that
$$ N_{Y,S}(P) %\le N_{\pi^{-1}(Y), S}(P')  
\le  \varepsilon C'	 h_{A}(P) + O(1)  $$
for all $P\in R \setminus Z$. Finally, one may eliminate the $O(1)$ (after possibly enlarging $Z$) by an application of Northcott's theorem for big divisors.
\end{proof}

Under an ampleness assumption, we repeatedly apply Theorem \ref{thm: count} to derive the stronger conclusion on the exceptional set in Theorem \ref{thm: exceptional}.

\begin{proof}[Proof of Theorem \ref{thm: exceptional}]
Let $Z$ be the Zariski closure of the set of points $P\in R$ such that
\begin{align*}
N_{Y,S}(P) \geq \varepsilon h_A(P).
\end{align*}
Since all of the assumptions of Theorem \ref{thm: count} are verified, $Z$ is a proper Zariski closed subset of $X$. After replacing $k$ by a finite extension, we may assume that every irreducible component of $Z$ is geometrically irreducible. Let $Z'$ be such an irreducible component and suppose that $Y|_{Z'}$ has codimension $\ge 2$ in $Z'$.

Let $\iota: Z' 	\hookrightarrow X$ be the embedding of $Z'$ in $X$ (identifying $Z'$ with its image where convenient). 
%Suppose $\dim Z' = m$. 
%We claim that there are at least $m + 1$ divisors among
%$D_1, \dots, D_{n+1}$ whose restriction to $Z'$ is 
Since $R$ is disjoint from the support of $D$, we may assume that $Z'$ is not contained in the support of any of the divisors $D_i$, and 
we can pull back the divisors $D_1,\dots, D_{n+1}$ to $Z'$. The numerical equivalences $d_jD_i|_{Z'} \equiv d_iD_j|_{Z'}$ are preserved and the divisors $D_i|_{Z'}$ remain ample for $i=1,\ldots, n+1$.  By Theorem \ref{thm: functoriality}, 
the set $R\cap Z'(k)$ is a $(D|_Z,S)$-integral set of points. Set $n' = \dim Z'$. 
Since $\bigcap_{i=1}^{n+1}D_i = \emptyset$, 
by Lemma 4.6 of \cite{Lev16}, 
the set of divisors $\{D_1|_{Z'},\dots,D_{n+1}|_{Z'}\}\subseteq \mathrm{Div}(Z') $ generate a
subgroup of rank at least $n'+1$ in $\mathrm{Div}(Z') $. 
Rearranging, we may assume that ${D_1|_{Z'},\dots,D_{n'+1}|_{Z'}}$
are linearly independent in $\mathrm{Div}(Z')$. 
%Since $Y\cap Z'$ is of codimension $\ge 2$ in $Z'$, 
Then by Theorem \ref{thm: count} applied to the variety $Z'$, the divisors 
$D_1|_{Z'}, \dots, D_{n'+1}|_{Z'}$, 
and the subscheme $\iota^*Y=Y|_{Z'}$, 
for any $\varepsilon' > 0$, 
there exists a proper Zariski closed subset $Z''\subset Z'$ such that 
$ N_{\iota^*Y,S}(P)< \varepsilon' h_{\iota^*A}(P)$
for all $P\in R\cap (Z'\setminus Z'')(k)$.  By Theorem \ref{thm: functoriality}, 
\begin{equation}\label{eq: small}
  N_{Y,S}(\iota(P)) < \varepsilon h_{A}(\iota(P)) + C'
\end{equation}
for some constant $C'$ and all $P\in R\cap (Z'\setminus Z'')(k)$.  Since $A$ is ample, we may eliminate the term $C'$ by adding finitely many points to $Z''$,  and the resulting inequality contradicts the definition of $Z'$. Thus, $Y|_{Z'}$ can have codimension at most one in $Z'$.
\end{proof}

We now give an estimate for the proximity function.  Compared with earlier results (e.g., \cite[Th.~3.1]{WY19}, \cite[Cor.~1.5]{HL20}) we do not require the divisors $D_i$ to be ample, nor do we require them to be in general position.

\begin{theorem}\label{thm: proximity}
In addition to the general setup, suppose that the closed subscheme $Y$ does not contain any point of the set
$\bigcup_{i=1}^{n+1}\bigcap_{j\neq i}\mathrm{Supp}D_j.$ 
Let $\varepsilon > 0$ and let $A$ be a big divisor on $X$. Then there exists a proper Zariski closed subset $Z$ of $X$ such that
\begin{equation} 
m_{Y,S}(P) < \varepsilon  h_A(P)
\end{equation}
for all $P \in R \setminus Z$. 
\end{theorem}

\begin{proof}%[Proof of Theorem \ref{thm: proximity}]
We may assume that $R$ is Zariski dense in $X$ (otherwise the conclusion of the theorem is trivial). Then by Lemma \ref{lem: reduction}, after replacing the divisors $D_i$ by suitable integer multiples, we assume that the divisors $D_i$ are linearly equivalent.  We define maps $\psi,\phi,\pi$ as in the proof of Theorem \ref{thm: count}, so that $\psi\circ \pi=\phi$ and the rational map $\psi:X\dashrightarrow \mathbb{P}^n$ is given by
$$\psi(P) = [f_1(P): \cdots: f_n(P):f_{n+1}(P)],$$
where $\mathrm{div}(f_i) = D_i -  D_{n+1}$, $i=1,\ldots, n+1$, and $f_{n+1}=1$. Noting that $\psi_i:=[f_1/f_i: \cdots: f_{n+1}/f_i]$ is defined for points outside $D_i$, $i=1,\ldots, n+1$, we see that $\psi$ may be extended to a morphism on $X\setminus \bigcap_{i=1}^{n+1} D_i$.

Now suppose $Y$ verifies the assumptions of the theorem. Since $Y$ does not pass through any point in $\bigcap_{i=1}^{n+1}D_i$, we have $\psi(Y) = \phi(\pi^{-1}(Y))$, which is closed as $\phi$ is proper.  Let $P_0  = [1 : 0 : \dots : 0],\dots, P_n = [0: 0 :\dots : 0 :1]$ in $\mathbb{P}^n(k)$.   If $P\not\in \Supp D_i$ and $\psi_i(P)$ has only one nonzero coordinate, then clearly $P\in \bigcap_{j\neq i}D_j$.  It follows that the image $\psi(Y)$ does not contain any of the points $P_i$, $i=0,\dots, n$. 

Then we may find an effective divisor $E$ on $\mathbb{P}^n$ such that $\psi(Y)\subseteq \Supp E$ and $P_0,\ldots, P_{n}\not\in \Supp E$. By Theorem 4.1 of \cite{Lev19}, 
there exists a finite union $Z'$ of translates of proper algebraic subgroups
of $\Gm^n$ such that
$m_{E,S}(P') \le \varepsilon \cdot h(P')$
for all points $P' \in \Gm^n(\mathcal{O}_{k,S}) \setminus Z' \subseteq \Pn(k)$. This implies that $m_{\pi^*Y,S}(P)\ll m_{\phi^*E,S}(P)\leq \varepsilon h_{\phi^*H}(P)$ for all $P\in \pi^{-1}(R)$ outside some proper Zariski closed subset of $X'$, where $H$ is a hyperplane in $\mathbb{P}^n$.  As in the proof of Theorem \ref{thm: count}, this in turn implies that there exists a proper Zariski closed subset $Z\subseteq X$ such that
$$m_{Y,S}(P)<\varepsilon h_A(P)$$ for all $P\in R\setminus Z$.
\end{proof}

\begin{proof}[Proof of Theorem \ref{thm: height}]
This is a straightforward combination of Theorem \ref{thm: count} and Theorem \ref{thm: proximity}. 
\end{proof}

We end this section by describing a family of surfaces $X\setminus D$ where: Theorem \ref{thm: count} applies,
 the surfaces $X\setminus D$ contain dense sets of integral points (so that the inequalities are not vacuous), and previous results do not na\"ively apply (of course, our results are ultimately derived from \cite{Lev19}).

\begin{example}
Consider the $n$th Hirzebruch surface $X := \mathbb{F}_n$, $n\geq 1$. 
%isomorphic to $ \mathbb{P}(\mathcal{O}_{\mathbb{P}^1} \oplus \mathcal{O}_{\mathbb{P}^1}(-2))$. 
Let $C_0$ denote a section of $\pi: X \rightarrow \mathbb{P}^1$, and let $f$ be a fiber of $\pi$. 
%Then as is well-known, $C_0^2 = -n, f^2 = 0$ and $(C_0.f) = 1$. 
We say that a divisor $D$ on $X$ is of type $(a,b)$ if  $D\sim aC_0 + bf$. 
Let $b\geq n$ be an integer and let $B_1,B_2,B_3,B_4$ be distinct irreducible curves of types $(0,1),(0,1),(1,b)$, and $(1,b+1)$, respectively, all defined over a number field $L$.  
Then we have three linearly independent (in $\mathrm{Div}(X)$) divisors
\begin{equation}
    \begin{aligned}
      D_1 & := B_1 + B_3, \\
      D_2 & := B_2 + B_3, \\
      D_3 & := B_4, 
    \end{aligned}
\end{equation}
all of type $(1,b+1)$. 
Let $D := D_1 + D_2 +D_3$, and let $B := B_1+B_2+B_3 + B_4$, a reduced divisor of type $(2,2b+3)$. By Theorem 5.5.4 of \cite{Lev05}, there exist choices of the curves $B_1,B_2,B_3,B_4$ such that $B$ is a normal crossings divisor and there exist dense sets of $D$-integral points over sufficiently large number
fields $M\supset L$. Note that the canonical divisor $K_X$ is of type $ (-2,-n-2)$, and in this case (when $B$ is a normal crossings divisor) the log-canonical divisor $K_X + B$ is of type $(0,2b-n+1)$. In particular,  since $2b-n+1\geq 1$, $X\setminus D$ has positive log Kodaira dimension, and $X\setminus D\not\cong \mathbb{G}_m^2$. Note also that the divisors $D_1,D_2,D_3$ are not in general position.  Let $S$ be a finite set of places of $L$ containing all the archimedean places, let $R \subset {X}(L)$ be a set of $(D,S)$-integral points, and let $A$ be a big divisor on $X$. 
Let $Y \subseteq {X}$ be a closed subscheme of codimension $2$. 
By Theorem \ref{thm: count}, 
there exists a proper Zariski closed subset $Z \subseteq {X}$
such that for all $P\in  R\setminus Z$
we have $N_{Y,S}(P) \le \varepsilon h_A(P)$.
 
As an explicit example, consider $\mathbb{F}_1$, which may be identified with the blowup $\pi:X\to \mathbb{P}^2$ at (say) $P=[0:0:1]\in \mathbb{P}^2$. Let $B_1$ and $B_2$ be the strict transforms (in $X$) of the lines defined by $x=0$ and $y=0$, respectively, let $B_3$ be the pullback of the line defined by $z=0$, and let $B_4$ be the strict transform of the conic defined by $y^2=x^2+xz$. Then these divisors fit into the framework above (with $b=1$). Let
\begin{align*}
R'=\left\{[uy:y:1]\mid u,v\in \mathcal{O}_{k,S}^*, y=\frac{u(v-1)}{u^2-1}\in \mathcal{O}_{k,S}\setminus\{0\}\right\}.
\end{align*}
Then $R=\pi^{-1}(R')\subset X(k)$ is a set of $(B_1+B_2+B_3+B_4,S)$-integral points since for $[x:y:1]\in R'$, we have that $x$ and $y$ are $S$-integers, and $x/y=u$ and $(y^2-x^2-x)/x=-v$ are $S$-units. We note that $R$ is Zariski dense in $X$; for instance, this follows already from considering the points in $R'$ where $v=u^{2n}$ for some even integer $n$ (so that $\frac{u(v-1)}{u^2-1}\in \mathcal{O}_{k,S}$ for all such $v$). Now Theorem \ref{thm: count} yields appropriate upper bounds on the gcd of two coprime polynomials evaluated at points of $R'$ (outside of some exceptional set).

%Now since
%\begin{equation}
    %\begin{aligned}
      %((1,3))^2 & = 2\cdot 3 - 2 = 4, 
    %\end{aligned}
%\end{equation}
%the map $\psi: X\setminus D  \rightarrow  \Gm^2$ is of degree $4$, 
%and hence we have gone beyond Example \ref{ex: unit}. 
%Hence $X\setminus D$ is not isomorphic to $\Gm^2$. 
\end{example}

\section{Some Unlikely Intersection Conjectures and Their Relation to GCD Inequalities}\label{unlikely}

We first state a general version of our main conjecture, which we will show follows from other standard conjectures (Theorem \ref{thm: unlikelyexceptional}).

\begin{conjecture}\label{conj: unlikely2}
Let $X$ be a nonsingular complex projective variety and let $D$ be a normal crossings divisor on $X$. Suppose that $\bar{\kappa}(X\setminus D)\geq 0$. Let $Y$ be a positive-dimensional subvariety of $X$ of codimension at least $2$ such that $Y$ intersects the support of $D$ transversally. Let $\mathbf{W}$ be the union of the subvarieties $W$ of $X\setminus D$ such that
\begin{enumerate}
\item $W$ is not of log general type. \label{conjcond1}
\item $Y\subset \overline{W}$ and $Y$ is an ample Cartier divisor on $\overline{W}$.\label{conjcond2}
\end{enumerate}
Then $\mathbf{W}$ is not Zariski dense in $X$.
\end{conjecture}

\begin{numberedremark}
\label{counterremark}
It is natural to raise the question of whether \eqref{conjcond2} may be replaced more generally by the condition that $Y\cap \overline{W}$ is positive-dimensional and an ample Cartier divisor on $\overline{W}$. We sketch a counterexample to this stronger formulation: if $S$ is an elliptic K3 surface with a section $\sigma_P$ of infinite order, then take $X=S\times S$, $Y=\Delta_S$ the diagonal, and $W_n=\sigma_{nP}\times \sigma_{nP}\cong \mathbb{P}^1\times\mathbb{P}^1$ where $\sigma_{nP}$ is a multiple of the section $\sigma_P$. In this case, $W_n\cap Y$ is the diagonal in $W_n$ (which is ample on $W_n\cong \mathbb{P}^1\times\mathbb{P}^1$) and the union of the $W_n$ is easily shown to be Zariski dense in $X$. We leave it as an open problem to find a more general formulation of Conjecture \ref{conj: unlikely2} in this direction.
\end{numberedremark}

\subsection{Conjecture \ref{conj: unlikely2} and Exceptional Sets in GCD Inequalities}

The goal of this section is to prove that, under conjectures of Lang and Griffiths (or Vojta), the set $\mathbf{W}$ of Conjecture \ref{conj: unlikely2} is contained in the exceptional set of a suitable ``GCD inequality"
\begin{align*}
T_{Y,f}(r)\leq_{\exc} \varepsilon T_{A,f}(r).
\end{align*}
In particular, Conjecture \ref{conj: unlikely2} is implied by other standard conjectures.

For generality, we work over $\mathbb{C}$ and in the setting of value distribution theory (the na\"ive analogous results using Diophantine conjectures would yield results only over $\overline{\mathbb{Q}}$).

We first state the two relevant conjectures.

We use the following version of a conjecture of Griffiths (analogous to Vojta's conjecture in Diophantine approximation):

\begin{conjecture}
\label{Vojtaconj}
Let $X$ be a nonsingular complex projective variety, $D$ a normal crossings divisor on $X$, $K$ the canonical divisor on $X$, $A$ a big divisor on $X$, and $\varepsilon>0$. There exists a proper Zariski closed subset $Z=Z(X,D,A,\varepsilon)$ of $X$ such that for all $C\in\mathbb{R}$, the inequality
\begin{align*}
m_{D,f}(r)+T_{K,f}(r)\leq_{\exc} \varepsilon T_{A,f}(r)+C
\end{align*}
holds for all nonconstant holomorphic maps $f:\mathbb{C}\to X$ whose image is not contained in $Z$.
\end{conjecture}

Recall that the notation $\leq_{\exc}$ means that the inequality holds for all $r>0$ outside a set of finite Lebesgue measure.

Before stating the second conjecture, we introduce some more terminology.

\begin{definition}
Let $X$ be a complex variety. We define the holomorphic exceptional set $\Exch(X)$ to be the Zariski closure of the union of all images of nonconstant holomorphic maps $f:\mathbb{C}\to X$.
\end{definition}

We have the following conjectural characterization of when the holomorphic exceptional set is all of $X$, giving a converse to the conjecture of Green-Griffiths (Lang \cite[Ch. I 3.5]{Lang91} conjectured this in the projective case).

\begin{conjecture}
\label{Langconj}
Let $X$ be a complex variety. Then $\Exch(X)$ is a proper Zariski closed subset of $X$ if and only if $X$ is of log general type.
\end{conjecture}

We now prove a general (conjectural) GCD inequality as a consequence of Conjecture \ref{Vojtaconj}.

\begin{theorem}
\label{Vgcd}
Assume Conjecture \ref{Vojtaconj}. Let $X$ be a nonsingular complex projective variety and let $D$ be a normal crossings divisor on $X$. Suppose that $\bar{\kappa}(X\setminus D)\geq 0$. Let $Y$ be a closed subscheme of $X$ of codimension at least $2$ such that the support of $Y$ intersects the support of $D$ transversally. Let $A$ be an ample divisor on $X$ and $\varepsilon>0$. Then there exists a proper Zariski closed subset $Z=Z(X,A,D,Y,\varepsilon)$ of $X$ such that the inequality
\begin{align*}
T_{Y,f}(r)\leq_{\exc} \varepsilon T_{A,f}(r)+O(1)
\end{align*}
holds for all nonconstant holomorphic maps $f:\mathbb{C}\to X\setminus D$ whose image is not contained in $Z$.
\end{theorem}

\begin{proof}
One easily reduces to the case that $Y$ is a subvariety of $X$. First suppose that $Y$ is nonsingular and of codimension $c\geq 2$ in $X$. Let $\pi:\tilde{X}\to X$ be the blowup along $Y$. Then we have the relation between canonical divisors \cite[Exercise II.8.5]{Har77},
\begin{align*}
K_{\tilde{X}}\sim \pi^*K_X+(c-1)E,
\end{align*}
where $E$ is the exceptional divisor on $\tilde{X}$ associated to the blowup. Since $Y$ intersects the support of $D$ transversally, $\pi^*D$ is again a normal crossings divisor. Let $f:\mathbb{C}\to X\setminus D$ be a nonconstant holomorphic map. Then $m_{D,f}(r)=T_{D,f}(r)+O(1)$. We may also assume $f(\mathbb{C})\not\subset Y$. Let $\tilde{f}:\mathbb{C}\to \tilde{X}$ be the lift of $f$ to the blowup. By Conjecture \ref{Vojtaconj}, there exists a proper Zariski closed subset $\tilde{Z}$ of $\tilde{X}$ such that
\begin{align*}
m_{\pi^*D,\tilde{f}}(r)+T_{K_{\tilde{X}},\tilde{f}}(r)&=T_{\pi^*D,\tilde{f}}(r)+T_{\pi^*K_X+(c-1)E,\tilde{f}}(r)+O(1)\\
&=T_{K_X+D,f}(r)+(c-1)T_{Y,f}(r)+O(1)\\
&\leq_{\exc} \varepsilon T_{\pi^*A,\tilde{f}}(r)+O(1)\\
&\leq_{\exc} \varepsilon T_{A,f}(r)+O(1)
\end{align*}
holds for all nonconstant holomorphic maps $f:\mathbb{C}\to X$ whose image is not contained in $Z=\pi(\tilde{Z})$. Since $\bar{\kappa}(X\setminus D)\geq 0$, after possibly enlarging $Z$, $T_{K_X+D,f}(r)$ is bounded from below, and the desired inequality follows.

The general case follows easily from this case using the functorial and additive properties of the characteristic function, along with embedded resolution of singularities by repeated blowups along smooth centers.
\end{proof}

Next, we prove a key technical result showing, roughly, that if a positive-dimensional subvariety $Y$ of $X$ is ample in some subvariety $Z$ of $X$, then $Y$ cannot be ``too small" in $Z$ (note also that the theorem fails badly when $Y$ is a point; see the related Example \ref{0dimexample}).
\begin{theorem}
\label{comparisontheorem}
Let $X$ be a nonsingular complex projective variety and let $A$ be an ample divisor on $X$.  Let $Y$ be a positive-dimensional subvariety of $X$.  There exists a constant $\delta>0$, depending only on $X$, $Y$, and $A$ with the following property: if $Z$ is a subvariety of $X$ such that $Y\subset Z$ and $Y$ is a nef and big Cartier divisor on $Z$, then there exists a proper Zariski closed subset $W$ of $Z$ and a constant $C\in \mathbb{R}$ such that for every holomorphic map $f:\mathbb{C}\to Z$ with image not contained in $W$,
\begin{align*}
T_{Y,f}(r)> \delta T_{A,f}(r)+C.
\end{align*}

\end{theorem}

\begin{proof}
Let $Z$ be a subvariety of $X$ of dimension $m$ such that $Y\subset Z$ and $Y$ is a nef and big Cartier divisor on $Z$. For clarity, we let $Y_Z$ denote this divisor on $Z$. Let $A_Z$ be an ample divisor on $Z$ such that $\mathcal{O}_Z(A_Z)=\mathcal{O}_X(A)|_Z$.

Since $Y\subset Z$ (see \cite[Prop.4]{Kle66} for a more general statement),
\begin{align*}
\left( A_Z^ {m-1} \cdot Y_Z \right) = (A^{m-1} \cdot Y),
\end{align*}
where the intersections take place on $Z$ and $X$, respectively.

By the generalized Hodge index theorem \cite[Corollary 1.6.3 (ii)]{Laz}, we have
\begin{align*}
(Y_Z^{m-1}\cdot A_Z) (A_Z^m)^{m-2}\leq (Y_Z\cdot A_Z^{m-1})^{m-1}=(Y\cdot A^{m-1})^{m-1}.
\end{align*}

Since $(A_Z^m)\geq 1$, we have
\begin{align*}
(Y_Z^{m-1}\cdot A_Z)\leq (Y\cdot A^{m-1})^{m-1}.
\end{align*}

Then for $N> m(Y\cdot A^{m-1})^{m-1}$,
\begin{align*}
((NY_Z)^m)&=N^m(Y_Z^m)> mN^{m-1}(Y\cdot A^{m-1})^{m-1}\\
&> m((NY_Z)^{m-1}\cdot A_Z).
\end{align*}

It follows from a criterion of Siu \cite[Th.~2.2.15]{Laz} that $NY_Z-A_Z$ is big. In particular, for some $M>0$, $M(NY_Z-A_Z)\sim E$ for some effective divisor $E$ on $Z$. Therefore, as long as the image of $f$ is not contained in the support of $E$, $T_{M(NY_Z-A_Z),f}(r)$ is bounded from below.
Since
\begin{align*}
T_{M(NY_Z-A_Z),f}(r)=M(NT_{Y_Z,f}(r)-T_{A_Z,f}(r))+O(1),
\end{align*}
we find that if the image of $f$ is not contained in the support of $E$, then for some constant $C$,
\begin{align*}
T_{Y_Z,f}(r)\geq \frac{1}{N}T_{A_Z,f}(r)+C.
\end{align*}
Thus, we may take $0<\delta<\frac{1}{N}$. Finally, we note that $N$ (and hence $\delta$) can be chosen depending only on $A$ and $Y$ (and not $Z$).
\end{proof}

\begin{theorem}
\label{thm: unlikelyexceptional}
Assume Conjecture \ref{Vojtaconj} and Conjecture \ref{Langconj}. Let $X$ be a nonsingular complex projective variety and let $D$ be a normal crossings divisor on $X$. Suppose that $\bar{\kappa}(X\setminus D)\geq 0$. Let $Y$ be a positive-dimensional subvariety of $X$ of codimension $\geq 2$ such that $Y$ intersects the support of $D$ transversally. 
Let $\mathbf{W}$ be the union of the subvarieties $W$ of $X\setminus D$ such that
\begin{enumerate}
\item $W$ is not of log general type. 
\item $Y\subset \overline{W}$ and $Y$ is a nef and big Cartier divisor on $\overline{W}$.
\end{enumerate}
Let $Z(X,A,D,Y,\varepsilon)$ be the proper Zariski closed subset of $X$ in Theorem \ref{Vgcd}. Then there exists $\delta>0$ such that $\mathbf{W}\subset Z(X,A,D,Y,\delta)$. In particular, $\mathbf{W}$ is not Zariski dense in $X$.
\end{theorem}

\begin{proof}
Let $\delta>0$ be as in Theorem \ref{comparisontheorem} (for $X$, $Y$, and $A$). Let $W$ be a subvariety of $X\setminus D$ satisfying \eqref{conjcond1} and \eqref{conjcond2}. By Theorem \ref{comparisontheorem}, there exists a proper Zariski closed subset $Z'$ of $W$ and a constant $C$ such that for every nonconstant holomorphic map $f:\mathbb{C}\to W$ with image not contained in $Z'$, 
\begin{align*}
T_{Y,f}(r)>\delta T_{A,f}(r)+C.
\end{align*}
Let $f:\mathbb{C}\to W$ be such a holomorphic map. Then since $A$ is ample and $f$ is nonconstant, for any constant $C'\in\mathbb{R}$, the inequality 
\begin{align*}
T_{Y,f}(r)\leq_{\exc} \frac{\delta}{2} T_{A,f}(r)+C'
\end{align*}
does not hold. Thus, if $Z=Z(X,A,D,Y,\frac{\delta}{2})$ is the proper Zariski closed subset of $X$ of Theorem \ref{Vgcd}, we have $f(\mathbb{C})\subset Z$. Since $W$ is not of log general type, as we vary $f$, the union of such images is Zariski dense in $W$ by Conjecture \ref{Langconj}. Thus, $W\subset Z$. It follows that $\mathbf{W}\subset Z=Z(X,A,D,Y,\frac{\delta}{2})$.
\end{proof}

\begin{corollary}
\label{toriccorollary}
Let $X$ be a smooth toric compactification of $\mathbb{G}_m^n$, so that $\mathbb{G}_m^n=X\setminus D$ for some effective divisor $D$ on the smooth projective toric variety $X$. Let $Y$ be a positive-dimensional subvariety of $X$ of codimension at least $2$ in $X$ such that $Y$ does not contain any point of intersection of $n$ distinct irreducible components of the boundary $X\setminus\mathbb{G}_m^n$.  Then Conjecture \ref{conj: unlikely2} holds (for $X$, $D$, and $Y$).
\end{corollary}

\begin{proof}
First, we note that Conjecture \ref{Langconj} holds for subvarieties of $\mathbb{G}_m^n$ (in fact, by Ueno-type results it holds more generally for subvarieties of semiabelian varieties). Moreover, by work of the second author \cite[Th.~1.14 (stated in the Diophantine setting)]{Lev19}, for $X$, $Y$ and $D$ as in the statement, Theorem \ref{Vgcd} is known unconditionally (i.e., without assuming Conjecture \ref{Vojtaconj}) under the weaker condition on $Y$ in the statement of the Corollary. Thus, in this case we may remove the use of Conjectures \ref{Vojtaconj} and \ref{Langconj} in the proof of Theorem \ref{thm: unlikelyexceptional} and obtain an unconditional proof of Conjecture \ref{conj: unlikely2} in this case.
\end{proof}

As a further illustration of Corollary \ref{toriccorollary}, we note the case of projective space.

\begin{corollary}
Let $D$ be the sum of the $n+1$ coordinate hyperplanes in $\mathbb{P}^n$. Let $Y$ be a positive-dimensional subvariety of $\mathbb{P}^n$ of codimension at least $2$ such that $Y$ does not contain any of the points $[1:0:\cdots :0], [0:1:0:\cdots:0],\ldots, [0:\cdots :0:1]$. 
Then Conjecture \ref{conj: unlikely2} holds (for $\mathbb{P}^n$, $D$, and $Y$).
\end{corollary}

\subsection{The Conjecture on Abelian Varieties}

%We prove Conjecture \ref{conj: unlikely3} for an abelian variety $X=A$ under the assumption that $3\dim Y<\dim A$.
%
%
%\begin{theorem}
%\label{abelian1}
%Let $A$ be a complex abelian variety and let $Y$ be a closed subvariety of $A$ of codimension $\ge 2$ satisfying 
%\begin{align*}
%3\dim Y<\dim A.
%\end{align*}
%Let $\mathbf{W}$ be the union of the subvarieties $W$ of $A$ such that
%\begin{enumerate}
%\item $W$ is not a variety of general type. 
%\item $Y\cap W$ is positive-dimensional and an ample Cartier divisor on $W$.
%\end{enumerate}
%Then $\mathbf{W}$ is not Zariski dense in $A$.
%\end{theorem}
%
%\begin{proof}
%Let $W$ satisfy \eqref{conjcond1} and \eqref{conjcond2}. Let 
%\begin{align*}
%B = \mathrm{stab}(W)= \{a\in A~|~a + W = W\} 
%\end{align*}
%be the stabilizer of $W$.
%Since $W$ is not of general type, by a result of Ueno \cite[Th.~4.2]{Voj96} we have that $\dim B > 0$. Consider the Ueno fibration $\pi:W \rightarrow W/B$. If $\pi(Y\cap W)\neq W/B$, then for any point $P\in W/B\setminus \pi(Y\cap W)$, the positive-dimensional fiber $\pi^{-1}(P)$ does not intersect $Y\cap W$, contradicting that $Y\cap W$ is ample in $W$. Therefore $W=B+(Y\cap W)$. Moreover, using ampleness again, for any $b\in B$, $(b+(Y\cap W))\cap (Y\cap W)\neq \emptyset$, and therefore $b\in (Y\cap W)-(Y\cap W)\subset Y-Y$. It follows that $W = B + (Y\cap W) \subset Y - Y + Y$, and $Y-Y+Y$ is a proper subvariety of $A$ of dimension at most $3 \dim Y < \dim A$. Therefore $\mathbf{W}$ is not Zariski dense in $A$.  
%\end{proof}

We prove Conjecture \ref{conj: unlikely2} when $X=A$ is an abelian variety.

\begin{theorem}
\label{abelian2}
Let $A$ be a complex abelian variety of dimension $g$, $D$ a very ample divisor on $A$, and $Y\subset A$ a subvariety of dimension $k\geq 1$ and codimension at least $2$ in $A$. Then there are only finitely many subvarieties $W$ of $A$ such that:
\begin{enumerate}
\item $W$ is not a variety of general type. 
\item $Y\subset W$ and $Y$ is an ample Cartier divisor on $W$.
\end{enumerate}
More precisely, let $d_Y=\deg_D Y=(Y\cdot D^k)$. Then every such $W$ is either the unique translate of an abelian subvariety containing $Y$ as an ample divisor, or $W=B+Y$ with $B={\rm Stab}(W)^0$, $0<\dim B\leq k$, and
\begin{align*}
\deg_DB\leq kd_Y^{2k+1}.
\end{align*}
In particular, by a result of R\'emond \cite[Prop.~4.1]{Remond00}, the number of such $W$ is bounded by $c(g, \deg_DA)d_Y^{16g^4}$ for an explicit constant $c(g, \deg_DA)$.
\end{theorem}

We note the stronger finiteness conclusion of the theorem (as compared to Conjecture \ref{conj: unlikely2}). Stronger, but more complicated, explicit bounds in the theorem can be obtained from the proof and  R\'emond \cite[Prop.~4.1]{Remond00}.

In the proof we use the following result bounding the degree of a generically finite map to a general type subvariety of an abelian variety; for $n=1$ the inequality is equivalent to Riemann-Hurwitz combined with the classical degree bound on the geometric genus $g(T)\leq \frac{(d-1)(d-2)}{2}$. The proof is a relatively easy combination of known results, which we give for lack of a suitable reference.

%For a smooth complex projective variety $X$ of dimension $n$, recall that the volume $\vol(K_X)=\limsup\limits_{m\to \infty} n!\,h^0(mK_X)/m^n$ is a birational invariant.

\begin{theorem}\label{Castelgen}
Let $T\subset\mathbb{P}^N$ be a projective variety of dimension $n\geq 1$ and degree $d$, let $U$ be a subvariety of an abelian variety $A$ with $U$ of general type and dimension $n$, and let $f:T\dashrightarrow U$ be a dominant rational map.  Then $d\ge n+3$, $T$ is of general type, and
\begin{equation}\label{eq:main}
  \deg f\ \le\ \frac{d\,(d-n-2)^{n}}{2\cdot n!}\ .
\end{equation}
\end{theorem}

\begin{proof}
Let $\widetilde T\to T$ and $\widetilde U\to U$ be resolutions of singularities.  Since $\dim T=\dim U$, we have that $f$ is generically finite and by \cite[Prop.~3.2]{GuerraPirola},

\begin{equation}\label{eq:GP}
  \deg(f)\,\vol(K_{\widetilde U})\ \le\ \vol(K_{\widetilde T}).
\end{equation}

Since $\widetilde U\to U\subset A$ is generically finite, $\widetilde U$ has maximal Albanese dimension, and by hypothesis $\widetilde U$ is of general type. Then by a result of Jiang \cite[Thm.~1.4]{Jiang},

\begin{equation}\label{eq:Jiang}
  \vol(K_{\widetilde U})\ \ge\ 2\cdot n!.
\end{equation}

We now compute an upper bound for $\vol(K_{\widetilde T})$. A general linear projection (or the identity if $T$ is already a hypersurface) gives a finite birational morphism $T\to X\subset\mathbb{P}^{n+1}$, where $X$ is a hypersurface of degree $d$, and the composite map $\pi:\widetilde T\to X$ is a desingularization. By singular adjunction for pluricanonical forms \cite[Cor.~3.7 and Lem.~4.2]{BeckSchicho}, for every $m\geq 1$,
\begin{align*}
\pi_*\omega_{\widetilde T}^{\otimes m} \hookrightarrow \mathcal{O}_X\left(m(d-n-2)\right).
\end{align*}
Therefore,
\begin{align*}
 h^0(\widetilde T,mK_{\widetilde T}) \leq h^0 \left(X,\mathcal{O}_X\left(m(d-n-2)\right)\right),
\end{align*}
and in particular $d\geq n+3$, since otherwise the right-hand side is bounded. This immediately implies the volume inequality
\begin{align}
\label{eq:cond}
\vol(K_{\widetilde T})\le d(d-n-2)^n,
\end{align}
Combining \eqref{eq:GP}--\eqref{eq:cond} yields the inequality.
\end{proof}

\begin{proof}[Proof of Theorem \ref{abelian2}]
Let $W$ satisfy \eqref{conjcond1} and \eqref{conjcond2}. 

We first claim that $Y$ is a variety of general type. If not, let $B'=\mathrm{stab}(Y)$. We have $\dim B'>0$ since $Y$ is not of general type. Consider the map $\pi':A\to A/B'$ restricted to $W$. We have $\dim \pi'(W)\geq \dim W-\dim B'$. On the other hand, since $Y$ is fibered by $B'$ over $\pi'(Y)$, we have $\dim \pi'(Y)=\dim Y-\dim B'<\dim \pi'(W)$ and $\dim \pi'(Y)<\dim Y$. This contraction of $Y$ contradicts the ampleness of $Y$ in $W$: if $\dim \pi'(W)<\dim W$, then a positive-dimensional fiber above $\pi'(W)\setminus \pi'(Y)$ does not intersect $Y$, and if $\dim \pi'(W)=\dim W$ then a general curve in $\pi'(W)$ does not intersect the codimension $2$ subset $\pi'(Y)$ (and the preimage of such a curve does not intersect $Y$). Thus, we may assume from now on that $Y$ is a variety of general type.

Let 
\begin{align*}
B = \mathrm{stab}(W)^0= \{a\in A~|~a + W = W\}^0,
\end{align*}
the identity component of the stabilizer of $W$. We note that there is at most one $W\supset Y$ such that $\dim B=\dim W$. If there were two distinct $W$'s (which must be translates of abelian subvarieties) then the intersection would be a translate of some abelian subvariety and equal to $Y$, contradicting that $Y$ is of general type. Thus, we may assume $\dim B<\dim W$.

Consider the morphism $\pi:W\to W/B$. We have $\dim W/B>0$ since $\dim B<\dim W$, and by a result of Ueno, $W/B$ is a variety of general type. Since $Y$ is ample in $W$, $\pi|_Y:Y\to W/B$ must be surjective (otherwise some positive-dimensional fiber of $\pi$ in $W$ does not intersect $Y$).  Let $r=\dim B\leq \dim Y$ and $n=\dim W/B$. We may assume $D$ intersects $W$ properly, and let $D_W=D|_W$. We view $Y$ as a divisor on $W$. Let $V$ be given as the intersection of $Y$ and $r-1$ general elements of $|D|$ (if $r=1$ then $V=Y$) so that $V$ is irreducible by Bertini's theorem. Let $F$ be a general fiber of the flat morphism $\pi$. Then $\deg \pi|_V=(F\cdot Y\cdot D_W^{r-1})$. By the proof of Theorem \ref{comparisontheorem}, taking $\delta=\frac{1}{(k+1)d_Y^k+1}$, for some effective $\mathbb{Q}$-divisor $E$, $Y\sim \delta D_W+E$. Then since $F$ is a fiber and $D$ is very ample,
\begin{align*}
\deg \pi|_V=(F\cdot Y\cdot D_W^{r-1})\geq \delta (F\cdot D_W^r)=\delta (F\cdot D^r)=\delta (B\cdot D^r)=\delta \deg_DB.
\end{align*}

On the other hand, since $\pi|_V:V\to W/B$ is dominant, by Theorem \ref{Castelgen}, $d_Y\geq n+3$ and
\begin{align*}
\deg \pi|_V\leq \frac{d_Y(d_Y-n-2)^n}{2n!}.
\end{align*}

Combining the inequalities, it follows that
\begin{align*}
\deg_D B\leq \frac{1}{\delta}\deg \pi|_V\leq \frac{d_Y(d_Y-n-2)^n((k+1)d_Y^k+1)}{2n!}.
\end{align*}
The last quantity is easily bounded by $kd_Y^{2k+1}$. By a result of Birkenhake and Lange \cite[Th.~4.1]{BL}, an abelian variety has only finitely many abelian subvarieties of bounded degree with respect to $D$, and an explicit bound was given by R\'emond \cite{Remond00}. Finally, we note that $B$ and $Y$ determine $W$ as $W=B+Y$ (since as noted earlier $\pi(Y)=W/B$). The last statement of the theorem follows from easy estimates using the bound in \cite[Prop.~4.1]{Remond00}.
\end{proof}

\subsection{Examples and Comments Related to the Conjectures}

We give several examples showing the necessity of the hypotheses in the conjectures, including in the weaker Conjecture \ref{conj: unlikely2}.

We first show that the condition requiring $X$ (or $X\setminus D$) to have nonnegative (log) Kodaira dimension is necessary. At the same time, nearly the same example shows that the union in Conjecture \ref{conj: unlikely2} may be Zariski dense if we remove the non-general type condition \eqref{conjcond1}.

\begin{example}
\label{KodairaExample}
Let $X=\bP^3$ and let $Y$ be a line in $X$. For any point $P\notin Y$, there exists a plane $H$ of $X$ that goes through $P$ and $Y$. Note that $H$ is not a variety of general type and $Y$ is ample in $H$. Hence the union $\mathbf{W}$ of all such planes $H$ is the whole space $X$. Note that $\kappa(X)=\kappa(\bP^3) = -\infty$, and so Conjecture \ref{conj: unlikely2} may fail for projective varieties of negative Kodaira dimension. 

As a slight variation on this example, let $D$ be the sum of the $4$ coordinate planes on $X=\mathbb{P}^3$ and consider $X\setminus D\cong \mathbb{G}_m^3$. Again, take $Y$ to be a line in $\mathbb{P}^3$, satisfying the additional condition that it intersects $D$ transversally. Then the union of the planes $H$ above will cover $X$, and $Y$ is ample in each $H$. This does not violate Conjecture \ref{conj: unlikely2} as $W=H\cap (X\setminus D)$ is of log general type for all but finitely many $H$ (for all but finitely many $H$, it is the complement of $4$ lines in general position in $H\cong \mathbb{P}^2$). This shows that one cannot drop condition \eqref{conjcond1} in Conjecture  \ref{conj: unlikely2}. 
\end{example}

The next example explains the requirement that $Y$ is positive-dimensional in Conjecture \ref{conj: unlikely2} (or Conjecture \ref{conj: unlikely3}).

\begin{example}
\label{0dimexample}
Let $E$ be an elliptic curve and let $X=E\times E$. Let $Y$ be the identity in $E\times E$. Then $Y$ is ample in every proper abelian subvariety of $X$ (since any nontrivial effective divisor on a curve is ample), and their union is Zariski dense in $X$.
\end{example}

We also show that in condition \eqref{conjcond2}, some kind of ampleness condition is necessary (though possibly ample may be replaced by big, or maybe 	nef and big).

\begin{example}
Let $X=E\times E\times E$, where $E$ is an elliptic curve. Let $Y$ be any curve in $X$. Let $E'$ be a one-dimensional abelian subvariety of $X$. Then for $E'\not\subset {\rm Stab}(Y)$, $Y$ is a divisor on $W_{E'}=E'+Y$ and $W_{E'}$ is not of general type. Since the union of such $E'$'s is Zariski dense in $X$, the union of the relevant $W_{E'}$'s is also Zariski dense in $X$. Thus, one cannot replace the condition \eqref{conjcond2} in Conjecture  \ref{conj: unlikely2} with just the condition that $Y$ is a divisor on $W$.
\end{example}

Lastly, we show that in the noncompact case, some condition on the intersection of $Y$ and $D$ is necessary. In particular, Conjecture \ref{conj: unlikely2} seems to depend delicately on the chosen compactification.
io
\begin{example}
\label{intersectionexample}
As in Example \ref{KodairaExample}, let $D$ be the sum of the $4$ coordinate planes on $X=\mathbb{P}^3$ and consider $X\setminus D\cong \mathbb{G}_m^3$. Let $Y\not\subset \Supp D$ be a line passing through the intersection of $3$ of the coordinate planes (so $Y$ does not intersect $D$ transversally). Let $H$ be a plane in $\mathbb{P}^3$ containing $Y$. Then $H\cap D$ consists of $3$ lines passing through a common point and a fourth line not containing the common point. It is easy to calculate that $H\setminus (H\cap D)\subset \mathbb{G}_m^3$ is not of log general type, $Y$ is ample in $H\cong \mathbb{P}^2$, and the union of such planes $H$ is all of $\mathbb{P}^3$. Therefore the conclusion of Conjecture \ref{conj: unlikely2} does not hold for $Y$.
\end{example}

\section{Proximity Functions on Toric Varieties with Applications}
\label{sec: toric}

In the final section, we give some refinements of Theorem \ref{thm: proximity}, studying the proximity function for toric varieties in our setting. We prove a result on the dimension of the exceptional set for the proximity function, first in the case of projective space, and then under our general setup. The proof for projective space exploits the toric structure on $\mathbb{P}^n$.  A \emph{toric variety} $X$ is a variety over $k$ containing $\mathbb{G}^n_m$ as a dense open subvariety such that the action of $\mathbb{G}^n_m$ on itself extends to an algebraic action of $\mathbb{G}^n_m$ on $X$ . 
We refer to \cite{CLS11} for more details on toric varieties. 

\begin{theorem}
\label{thm: toric}
Let $X$ be a projective toric variety of dimension $n$ defined over a number field $k$. Let $A$ be an ample divisor on $X$. Let $Y$ be a closed subscheme of $X$, defined over $k$, whose support does not contain a torus fixed point. Let $S$ be a finite set of places of $k$ including the archimedean ones and let $\varepsilon > 0$. Then there exists a finite union $Z$ of translates of proper algebraic subgroups of $\Gm^n$ such that
$$m_{Y,S}(P) \le \varepsilon h_A(P) + O(1)$$ for all points 
$ P \in \Gm^n(\mathcal{O}_{k,S}) \setminus Z \subset X(k)$.
\end{theorem}

\begin{proof}
We may find an effective Cartier divisor $D$ on $X$ such that $D$ does not contain any torus fixed point and $Y\subset D$ as closed subschemes. In this case, $m_{Y,S}(P)\leq m_{D,S}(P)+O(1)$, and we may reduce to the case that $Y=D$ is a divisor. By Proposition 3.9.1 of \cite{Kollar}, there exists an equivariant resolution of singularities $\pi: \tilde{X} \rightarrow X$.  Since $\pi$ is equivariant, 
the pullback $\pi^*D$ does not contain a torus fixed point either.
%Hence $\pi^*D$ is in general position with respect to the boundary of the pair $(X,X \setminus \Gm^n)$. 

Now we apply the proof of Theorem 4.4 of \cite{Lev19} to the divisor $\pi^*D$ and the big divisor $\pi^*A$ to obtain
\begin{equation}
    m_{\pi^*D,S}(P) \le \varepsilon h_{\pi^*A}(P) + O(1)
\end{equation}
for all points $P \in \Gm^n(\mathcal{O}_{k,S}) \setminus \tilde{Z} \subset \tilde{X}(k)$, where  $\tilde{Z}$ is a finite union of
translates of proper algebraic subgroups of $\Gm^n$. The desired inequality then follows from functoriality. 
\end{proof}

We give a more refined characterization of the exceptional set in the following theorem. %Note that the assumption that $\ell < n$ is ``generically" satisfied.  

\begin{theorem}\label{thm: excToric}
Let $X$ be a projective toric variety of dimension $n$ defined over a number field $k$. Let $A$ be an ample divisor on $X$. 
Let $Y$ be a closed subscheme of $X$ defined over $k$. 
%Let $D$ be an effective divisor on $X$, defined over $K$, that is in general position with the boundary of $\Gm^n$ in $X$.
Let $S$ be a finite set of places of $k$ containing the archimedean places and let $\varepsilon > 0$. 
Let $\ell$ be the maximum non-negative integer such that there exists an orbit $O(P), P\in X(\Kbar)$, of codimension $\ell$ that intersects $Y$. %, and assume that $\ell < n$. 
Then there exists a finite union $Z$ of translates of proper algebraic subgroups of $\Gm^n$
of dimension at most $\ell$
such that
$$m_{Y,S}(P) \le \varepsilon h_A(P) + O(1)$$ for all points 
$ P \in \Gm^n(\mathcal{O}_{k,S}) \setminus (Y\cup Z) \subset X(k)$.
\end{theorem}

\begin{proof}
If $\ell \ge n$, then the result is vacuous. 
So we may assume that $\ell < n$. This implies that $Y$ does not contain any torus fixed point of $X$.
By Theorem \ref{thm: toric}, 
there exists a finite union $Z'$ of translates of proper algebraic subgroups of $\Gm^n$
such that 
$$m_{Y,S}(P) \le \varepsilon h_A(P) + O(1)$$ 
holds for all points 
$ P \in \Gm^n(\mathcal{O}_{k,S}) \setminus Z' \subset X(k)$.

Let $W=gT$ be such a translate in $Z'$. After possibly replacing $k$ by a finite extension (and replacing $S$ by the set of places lying above places of $S$), we may assume that $T$ is geometrically irreducible.  Clearly, by induction, it suffices to show that if $\dim W>\ell$ and $W\not\subseteq Y$, then there exists a proper subset $Z''\subset W$ such that $Z''$ is a finite union of translates of proper algebraic subgroup of $\Gm^n$ and $$m_{Y,S}(P) \le \varepsilon h_A(P) + O(1)$$ 
for all points  $P \in W(k)\cap \Gm^n(\mathcal{O}_{k,S}) \setminus Z'' \subset X(k)$.

By enlarging $S$, we may assume that $g\in  \Gm^n(\mathcal{O}_{k,S})$. Then $gP\in (gT)(k)\cap  \Gm^n(\mathcal{O}_{k,S})$ if and only if $P\in T(k)\cap  \Gm^n(\mathcal{O}_{k,S})$.  We also have $$m_{Y,S}(gP)=m_{g^*Y,S}(P)+O(1), ~h_A(gP)=h_{g^*A}(P)+O(1).$$
An orbit $O(P)$ intersects $Y$ if and only if it intersects $g^*Y$.  Then after replacing $Y$ by $g^*Y$ and $W=gT$ by $T$, it suffices to prove the desired result when $W=T\cong \mathbb{G}_m^{\dim T}$ is a subtorus of $\mathbb{G}_m^n$. Let $\overline{T}$ be the Zariski closure of $T$ in $X$.  Then $T$ acts naturally on $\overline{T}$, and $\overline{T}$ has an induced structure of a toric variety. We claim that $Y|_{\overline{T}}$ does not contain any $T$-invariant point.  Indeed, suppose that $P\in \overline{T}(\Kbar)$ is such a $T$-invariant point. 
Consider the orbit $O(P)$ under the action of $\Gm^n$. Let $\mathrm{stab}(P)$ be the stabilizer of $P$. Since $\mathrm{stab}(P)\supset T$, we have that $\mathrm{stab}(P)$ has dimension at least $\dim T$. 
Hence the orbit $O(P)$ under the action of $\Gm^n$ has codimension (in $X$) at least $\dim T > \ell$, contradicting the definition of $\ell$.  Then functoriality and Theorem \ref{thm: toric} applied to $\overline{T}$ and $Y|_{\overline{T}}$ give the desired result.
\end{proof}

We have the following important immediate corollary for the case of $\mathbb{P}^n$. 

\begin{corollary}\label{thm: excPn}
Let $Y$ be a closed subscheme of $\mathbb{P}^n$ defined over a number field $k$. Let $H_1,\ldots, H_{n+1}$ be the $n+1$ coordinate hyperplanes in $\mathbb{P}^n$. Suppose that $\bigcap_{i\in I}H_i\bigcap Y=\emptyset$ for every subset $I\subset \{1,\ldots, n+1\}$ with $|I|>\ell$.  Let $S$ be a finite set of places of $k$ containing the archimedean places. 
For any $\varepsilon > 0$, there exists a finite union $Z$ of translates of proper algebraic subgroups of $\Gm^n$ of dimension at most $\ell$ such that
$$m_{Y,S}(P) \le \varepsilon h(P) + O(1) $$
for all points $P \in \Gm^n(\mathcal{O}_{k,S}) \setminus (Y\cup Z) \subset \mathbb{P}^n(k)$.
\hfill \qedsymbol
\end{corollary}

More generally, we consider the exceptional set for the proximity function in the general case.

\begin{theorem}\label{thm: proxExc}
In addition to the general setup, 
let $\varepsilon > 0$ and let $A$ be an ample divisor on $X$. 
Assume that $\bigcap_{i=1}^{n+1} D_i = \emptyset$ and that $D_1,\ldots, D_{n+1}$ are ample. Suppose that 
\begin{align}
\label{inthyp}
\bigcap_{i\in I}D_i\bigcap Y=\emptyset
\end{align}
for every subset $I\subset \{1,\ldots, n+1\}$ with $|I|>\ell$.   Then there exists a Zariski closed subset $Z$ of $X$ of dimension at most
$\ell$ such that
\begin{equation*} %\label{eq: proxExc}
m_{Y,S}(P) < \varepsilon  h_A(P)
\end{equation*}
for all $P \in R \setminus (Y\cup Z)$. 
\end{theorem}

\begin{proof}
Let $\overline{R}$ be the Zariski closure of $R$ in $X$. After possibly extending $k$, we may assume that every irreducible component of $\overline{R}$ (over $k$) is geometrically irreducible. After replacing $D_i$ by an appropriate multiple, we may further assume that $d_i=1$ for all $i$.  Let $X'\not\subseteq Y$ be an irreducible component of $\overline{R}$, and let $R'=R\cap X'$, $Y'=Y|_{X'}$, $D_i'=D_i|_{X'}$, $A'=A|_{X'}$. By functoriality, it suffices to show that there exists a proper Zariski closed subset $Z$ of $X'$ of dimension at most $\ell$ such that
\begin{equation*}% \label{eq: proxExc}
m_{Y',S}(P) < \varepsilon  h_{A'}(P)
\end{equation*}
for all $P \in R' \setminus Z$. Let $n'=\dim X'$. Since $\bigcap_{i=1}^{n+1}D_i = \emptyset$, 
it follows from \cite[Lemma 4.6]{Lev16} that $D_1',\ldots, D_{n+1}'$ generate a
subgroup of rank at least $n'+1$ in $\mathrm{Div}(X') $. After reindexing, suppose that $D_1',\ldots, D_{n'+1}'$ are $\mathbb{Z}$-linearly independent in $\mathrm{Div}(X')$.  Since $R'$ is a Zariski dense set of $(D',S)$-integral points on $X'$, where $D'=\sum_{i=1}^{n'+1}D_i'$, it follows from the proof of Lemma \ref{lem: reduction} that $q(X')=0$, $\Pic(X')\cong \NS(X')$ and (after replacing the divisors by suitable multiples) $D_i'\sim D_j'$ for all $i$ and $j$. As $\bigcap_{i=1}^{n+1}D_i' = \emptyset$ as well, we can construct a morphism $\phi:X'\to\mathbb{P}^{n}$ such that if $H_1,\ldots, H_{n+1}$ are the $n+1$ coordinate hyperplanes of $\mathbb{P}^{n}$, then $\phi^*H_i=D_i'$ for $i=1,\ldots, n+1$.  It follows that $\bigcap_{i\in I}H_i\bigcap \phi(Y')=\emptyset$ for every subset $I\subset \{1,\ldots, n+1\}$ with $|I|>\ell$. Furthermore, after possibly enlarging $S$, $\phi(R')\subset \mathbb{G}_m^{n}(\mathcal{O}_{k,S})$. Then by Corollary \ref{thm: excPn} applied to the closed subscheme $\phi(Y')$, Lemma \ref{lem: pullback2}, and Theorem \ref{thm: functoriality}, there exists a closed subset $W\subset  \mathbb{G}_m^{n}, \dim W\leq \ell$, such that
\begin{equation*}% \label{eq: proxExc}
m_{Y',S}(P) < \varepsilon  h_{A'}(P)
\end{equation*}
for all $P \in R' \setminus \phi^{-1}(W\cup \phi(Y'))$. Finally, since $D_1',\ldots, D_{n+1}'$ are ample, $\phi$ is a finite morphism \cite[Cor.~1.2.15]{Laz}. Noting that since the $D_i$ are ample, \eqref{inthyp} implies $\dim Y\leq \ell$ (and so $\dim \phi(Y')\leq \ell$), we find that $\dim \phi^{-1}(W\cup \phi(Y'))\leq \ell$, as desired.
\end{proof}

\section*{Acknowledgment}
We would like to thank Zheng Xiao for many helpful discussions. 
We also thank Julie Wang for useful discussions, especially on Nevanlinna theory. 
We thank Nathan Grieve for helpful suggestions on the references. The counterexample in Remark \ref{counterremark} was suggested independently by ChatGPT and Claude. The authors subsequently verified the construction and the argument. No generative-AI system was otherwise used in obtaining the mathematical results of this article.

\bibliography{refFile}{}

@incollection {Klei05,
    AUTHOR = {Kleiman, Steven L.},
     TITLE = {The {P}icard scheme},
 BOOKTITLE = {Fundamental algebraic geometry},
    SERIES = {Math. Surveys Monogr.},
    VOLUME = {123},
     PAGES = {235--321},
 PUBLISHER = {Amer. Math. Soc., Providence, RI},
      YEAR = {2005},
      ISBN = {0-8218-3541-6},
   MRCLASS = {14C22},
  MRNUMBER = {2223410},
}

@article {CZ05,
    AUTHOR = {Corvaja, Pietro and Zannier, Umberto},
     TITLE = {A lower bound for the height of a rational function at
              {$S$}-unit points},
   JOURNAL = {Monatsh. Math.},
  FJOURNAL = {Monatshefte f\"{u}r Mathematik},
    VOLUME = {144},
      YEAR = {2005},
    NUMBER = {3},
     PAGES = {203--224},
      ISSN = {0026-9255},
   MRCLASS = {11J25 (11G50)},
  MRNUMBER = {2130274},
MRREVIEWER = {Yann Bugeaud},
       DOI = {10.1007/s00605-004-0273-0},
       URL = {https://doi-org.proxy1.cl.msu.edu/10.1007/s00605-004-0273-0},
}

@article {BCZ,
    AUTHOR = {Bugeaud, Yann and Corvaja, Pietro and Zannier, Umberto},
     TITLE = {An upper bound for the {G}.{C}.{D}. of {$a^n-1$} and
              {$b^n-1$}},
   JOURNAL = {Math. Z.},
  FJOURNAL = {Mathematische Zeitschrift},
    VOLUME = {243},
      YEAR = {2003},
    NUMBER = {1},
     PAGES = {79--84},
      ISSN = {0025-5874},
   MRCLASS = {11J25 (11D75)},
  MRNUMBER = {1953049},
MRREVIEWER = {Damien Roy},
       DOI = {10.1007/s00209-002-0449-z},
       URL = {https://doi-org.proxy2.cl.msu.edu/10.1007/s00209-002-0449-z},
}

@book {Har77,
    AUTHOR = {Hartshorne, Robin},
     TITLE = {Algebraic geometry},
      NOTE = {Graduate Texts in Mathematics, No. 52},
 PUBLISHER = {Springer-Verlag, New York-Heidelberg},
      YEAR = {1977},
     PAGES = {xvi+496},
      ISBN = {0-387-90244-9},
   MRCLASS = {14-01},
  MRNUMBER = {0463157},
MRREVIEWER = {Robert Speiser},
}

@article{HL20,
  title={A generalized Schmidt subspace theorem for closed subschemes},
  author={Heier, Gordon and Levin, Aaron},
  journal={American Journal of Mathematics},
  volume={143},
  number={1},
  pages={213--226},
  year={2021},
  publisher={Johns Hopkins University Press}
}

@article {Lev05,
    AUTHOR = {Levin, Aaron},
     TITLE = {Generalizations of {S}iegel's and {P}icard's Theorems},
   JOURNAL = {PhD thesis, University of California, Berkeley},
  FJOURNAL = {},
    VOLUME = {},
      YEAR = {2005},
    NUMBER = {},
     PAGES = {},
      ISSN = {},
   MRCLASS = {},
  MRNUMBER = {},
MRREVIEWER = {},
       DOI = {},
       URL = {},
}

@article {Lev16,
    AUTHOR = {Levin, Aaron},
     TITLE = {Integral points of bounded degree on affine curves},
   JOURNAL = {Compos. Math.},
  FJOURNAL = {Compositio Mathematica},
    VOLUME = {152},
      YEAR = {2016},
    NUMBER = {4},
     PAGES = {754--768},
      ISSN = {0010-437X},
   MRCLASS = {11G30 (14H25 30D35)},
  MRNUMBER = {3484113},
MRREVIEWER = {Joseph H. Silverman},
       DOI = {10.1112/S0010437X15007708},
       URL = {https://doi-org.proxy1.cl.msu.edu/10.1112/S0010437X15007708},
}

@article {Lev19,
    AUTHOR = {Levin, Aaron},
     TITLE = {Greatest common divisors and {V}ojta's conjecture for blowups
              of algebraic tori},
   JOURNAL = {Invent. Math.},
  FJOURNAL = {Inventiones Mathematicae},
    VOLUME = {215},
      YEAR = {2019},
    NUMBER = {2},
     PAGES = {493--533},
      ISSN = {0020-9910},
   MRCLASS = {11J25 (11B37 11G35)},
  MRNUMBER = {3910069},
MRREVIEWER = {Sajad Salami},
       DOI = {10.1007/s00222-018-0831-z},
       URL = {https://doi-org.proxy1.cl.msu.edu/10.1007/s00222-018-0831-z},
}

@article {LW20,
    AUTHOR = {Levin, Aaron and Wang, Julie Tzu-Yueh},
     TITLE = {Greatest common divisors of analytic functions and
              {N}evanlinna theory on algebraic tori},
   JOURNAL = {J. Reine Angew. Math.},
  FJOURNAL = {Journal f\"{u}r die Reine und Angewandte Mathematik. [Crelle's
              Journal]},
    VOLUME = {767},
      YEAR = {2020},
     PAGES = {77--107},
      ISSN = {0075-4102},
   MRCLASS = {32H30 (11J97)},
  MRNUMBER = {4160303},
MRREVIEWER = {Tuan Dinh Huynh},
       DOI = {10.1515/crelle-2019-0033},
       URL = {https://doi-org.proxy2.cl.msu.edu/10.1515/crelle-2019-0033},
}

@article {Nog81,
    AUTHOR = {Noguchi, Junjiro},
     TITLE = {Lemma on logarithmic derivatives and holomorphic curves in
              algebraic varieties},
   JOURNAL = {Nagoya Math. J.},
  FJOURNAL = {Nagoya Mathematical Journal},
    VOLUME = {83},
      YEAR = {1981},
     PAGES = {213--233},
      ISSN = {0027-7630},
   MRCLASS = {32H25 (32D15)},
  MRNUMBER = {632655},
MRREVIEWER = {Pierre Lelong},
       URL = {http://projecteuclid.org.proxy2.cl.msu.edu/euclid.nmj/1118786486},
}

@article {RW20,
    AUTHOR = {Ru, Min and Vojta, Paul},
     TITLE = {A birational {N}evanlinna constant and its consequences},
   JOURNAL = {Amer. J. Math.},
  FJOURNAL = {American Journal of Mathematics},
    VOLUME = {142},
      YEAR = {2020},
    NUMBER = {3},
     PAGES = {957--991},
      ISSN = {0002-9327},
   MRCLASS = {14C20 (11G35 11J97 32H30)},
  MRNUMBER = {4101336},
       DOI = {10.1353/ajm.2020.0022},
       URL = {https://doi-org.proxy1.cl.msu.edu/10.1353/ajm.2020.0022},
}

@book {NW14,
    AUTHOR = {Noguchi, Junjiro and Winkelmann, J\"{o}rg},
     TITLE = {Nevanlinna theory in several complex variables and
              {D}iophantine approximation},
    SERIES = {Grundlehren der Mathematischen Wissenschaften [Fundamental
              Principles of Mathematical Sciences]},
    VOLUME = {350},
 PUBLISHER = {Springer, Tokyo},
      YEAR = {2014},
     PAGES = {xiv+416},
      ISBN = {978-4-431-54570-5; 978-4-431-54571-2},
   MRCLASS = {32-02 (11J25 11J97 32H30 32Q45)},
  MRNUMBER = {3156076},
MRREVIEWER = {William A. Cherry},
       DOI = {10.1007/978-4-431-54571-2},
       URL = {https://doi-org.proxy1.cl.msu.edu/10.1007/978-4-431-54571-2},
}

@article {NWY02,
    AUTHOR = {Noguchi, Junjiro and Winkelmann, J\"{o}rg and Yamanoi, Katsutoshi},
     TITLE = {The second main theorem for holomorphic curves into
              semi-abelian varieties},
   JOURNAL = {Acta Math.},
  FJOURNAL = {Acta Mathematica},
    VOLUME = {188},
      YEAR = {2002},
    NUMBER = {1},
     PAGES = {129--161},
      ISSN = {0001-5962},
   MRCLASS = {32H30},
  MRNUMBER = {1947460},
MRREVIEWER = {William A. Cherry},
       DOI = {10.1007/BF02392797},
       URL = {https://doi-org.proxy2.cl.msu.edu/10.1007/BF02392797},
}

@article {Sil87,
    AUTHOR = {Silverman, Joseph H.},
     TITLE = {Arithmetic distance functions and height functions in
              {D}iophantine geometry},
   JOURNAL = {Math. Ann.},
  FJOURNAL = {Mathematische Annalen},
    VOLUME = {279},
      YEAR = {1987},
    NUMBER = {2},
     PAGES = {193--216},
      ISSN = {0025-5831},
   MRCLASS = {11G99 (11D99 14G99)},
  MRNUMBER = {919501},
MRREVIEWER = {Andrew Bremner},
       DOI = {10.1007/BF01461718},
       URL = {https://doi-org.proxy1.cl.msu.edu/10.1007/BF01461718},
}

@article{Sil05,
    AUTHOR = {Silverman, Joseph H.},
     TITLE = {Generalized greatest common divisors, divisibility sequences,
              and {V}ojta's conjecture for blowups},
   JOURNAL = {Monatsh. Math.},
  FJOURNAL = {Monatshefte f\"ur Mathematik},
    VOLUME = {145},
      YEAR = {2005},
    NUMBER = {4},
     PAGES = {333--350},
      ISSN = {0026-9255},
     CODEN = {MNMTA2},
   MRCLASS = {11G35 (11D75 11J25 14G25)},
  MRNUMBER = {2162351},
MRREVIEWER = {Yann Bugeaud},
       DOI = {10.1007/s00605-005-0299-y},
       URL = {http://dx.doi.org/10.1007/s00605-005-0299-y},
}

@article {WY19,
    AUTHOR = {Wang, Julie Tzu-Yueh and Yasufuku, Yu},
     TITLE = {Greatest common divisors of integral points of numerically
              equivalent divisors},
   JOURNAL = {Algebra Number Theory},
  FJOURNAL = {Algebra \& Number Theory},
    VOLUME = {15},
      YEAR = {2021},
    NUMBER = {1},
     PAGES = {287--305},
      ISSN = {1937-0652},
   MRCLASS = {11J97 (11J87 14G05 32A22)},
  MRNUMBER = {4226990},
       DOI = {10.2140/ant.2021.15.287},
       URL = {https://doi-org.proxy1.cl.msu.edu/10.2140/ant.2021.15.287},
}

@book {Voj87,
    AUTHOR = {Vojta, Paul},
     TITLE = {Diophantine approximations and value distribution theory},
    SERIES = {Lecture Notes in Mathematics},
    VOLUME = {1239},
 PUBLISHER = {Springer-Verlag, Berlin},
      YEAR = {1987},
     PAGES = {x+132},
      ISBN = {3-540-17551-2},
   MRCLASS = {11G35 (11J68 14G25 30D35 32H30)},
  MRNUMBER = {883451},
MRREVIEWER = {J. Silverman},
       DOI = {10.1007/BFb0072989},
       URL = {https://doi-org.proxy1.cl.msu.edu/10.1007/BFb0072989},
}

@article {Voj96,
    AUTHOR = {Vojta, Paul},
     TITLE = {Integral points on subvarieties of semiabelian varieties. {I}},
   JOURNAL = {Invent. Math.},
  FJOURNAL = {Inventiones Mathematicae},
    VOLUME = {126},
      YEAR = {1996},
    NUMBER = {1},
     PAGES = {133--181},
      ISSN = {0020-9910},
   MRCLASS = {14G05 (11G10)},
  MRNUMBER = {1408559},
MRREVIEWER = {Alexandru Buium},
       DOI = {10.1007/s002220050092},
       URL = {https://doi-org.proxy1.cl.msu.edu/10.1007/s002220050092},
}

@book {CLS11,
    AUTHOR = {Cox, David A. and Little, John B. and Schenck, Henry K.},
     TITLE = {Toric varieties},
    SERIES = {Graduate Studies in Mathematics},
    VOLUME = {124},
 PUBLISHER = {American Mathematical Society, Providence, RI},
      YEAR = {2011},
     PAGES = {xxiv+841},
      ISBN = {978-0-8218-4819-7},
   MRCLASS = {14M25 (05A15 05E45 52B12)},
  MRNUMBER = {2810322},
MRREVIEWER = {Ivan Arzhantsev},
       DOI = {10.1090/gsm/124},
       URL = {https://doi-org.proxy2.cl.msu.edu/10.1090/gsm/124},
}

@book {Kollar,
    AUTHOR = {Koll\'{a}r, J\'{a}nos},
     TITLE = {Lectures on resolution of singularities},
    SERIES = {Annals of Mathematics Studies},
    VOLUME = {166},
 PUBLISHER = {Princeton University Press, Princeton, NJ},
      YEAR = {2007},
     PAGES = {vi+208},
      ISBN = {978-0-691-12923-5; 0-691-12923-1},
   MRCLASS = {14E15 (32S45)},
  MRNUMBER = {2289519},
MRREVIEWER = {Dan Abramovich},
}

@book {Laz,
    AUTHOR = {Lazarsfeld, Robert},
     TITLE = {Positivity in algebraic geometry. {I}},
    SERIES = {Ergebnisse der Mathematik und ihrer Grenzgebiete. 3. Folge. A
              Series of Modern Surveys in Mathematics [Results in
              Mathematics and Related Areas. 3rd Series. A Series of Modern
              Surveys in Mathematics]},
    VOLUME = {48},
      NOTE = {Classical setting: line bundles and linear series},
 PUBLISHER = {Springer-Verlag, Berlin},
      YEAR = {2004},
     PAGES = {xviii+387},
      ISBN = {3-540-22533-1},
   MRCLASS = {14-02 (14C20)},
  MRNUMBER = {2095471},
MRREVIEWER = {Mihnea Popa},
       DOI = {10.1007/978-3-642-18808-4},
       URL = {https://doi-org.proxy2.cl.msu.edu/10.1007/978-3-642-18808-4},
}

@book {GW,
    AUTHOR = {G\"{o}rtz, Ulrich and Wedhorn, Torsten},
     TITLE = {Algebraic geometry {I}. {S}chemes---with examples and
              exercises},
    SERIES = {Springer Studium Mathematik---Master},
      NOTE = {Second edition [of  2675155]},
 PUBLISHER = {Springer Spektrum, Wiesbaden},
      YEAR = {[2020] \copyright 2020},
     PAGES = {vii+625},
      ISBN = {978-3-658-30732-5; 978-3-658-30733-2},
   MRCLASS = {14-01},
  MRNUMBER = {4225278},
       DOI = {10.1007/978-3-658-30733-2},
       URL = {https://doi-org.proxy2.cl.msu.edu/10.1007/978-3-658-30733-2},
}

@article {CLZ22,
    AUTHOR = {Corvaja, Pietro and Levin, Aaron and Zannier, Umberto},
     TITLE = {Intersections in subvarieties of {${\Bbb{G}}^l_{\rm{m}}$} and
              applications to lacunary polynomials},
   JOURNAL = {Trans. Amer. Math. Soc.},
  FJOURNAL = {Transactions of the American Mathematical Society},
    VOLUME = {375},
      YEAR = {2022},
    NUMBER = {1},
     PAGES = {71--86},
      ISSN = {0002-9947},
   MRCLASS = {14M25 (11J25 14D07)},
  MRNUMBER = {4358662},
MRREVIEWER = {Zhizhong Huang},
       DOI = {10.1090/tran/8470},
       URL = {https://doi-org.proxy2.cl.msu.edu/10.1090/tran/8470},
}

@article{Kle66,
 author = {Kleiman, S. L.},
 title = {Toward a numerical theory of ampleness},
 fjournal = {Annals of Mathematics. Second Series},
 journal = {Ann. Math. (2)},
 issn = {0003-486X},
 volume = {84},
 pages = {293--344},
 year = {1966},
 language = {English},
 doi = {10.2307/1970447},
 zbMATH = {3235756},
 Zbl = {0146.17001}
}

@article {BMZ,
    AUTHOR = {Bombieri, E. and Masser, D. and Zannier, U.},
     TITLE = {Anomalous subvarieties---structure theorems and applications},
   JOURNAL = {Int. Math. Res. Not. IMRN},
  FJOURNAL = {International Mathematics Research Notices. IMRN},
      YEAR = {2007},
    NUMBER = {19},
     PAGES = {Art. ID rnm057, 33},
      ISSN = {1073-7928,1687-0247},
   MRCLASS = {11G10 (11G35 11G50 14L10)},
  MRNUMBER = {2359537},
MRREVIEWER = {Dimitrios\ Poulakis},
       DOI = {10.1093/imrn/rnm057},
       URL = {https://doi-org.proxy2.cl.msu.edu/10.1093/imrn/rnm057},
}

@book {Zan12,
    AUTHOR = {Zannier, Umberto},
     TITLE = {Some problems of unlikely intersections in arithmetic and
              geometry},
    SERIES = {Annals of Mathematics Studies},
    VOLUME = {181},
      NOTE = {With appendixes by David Masser},
 PUBLISHER = {Princeton University Press, Princeton, NJ},
      YEAR = {2012},
     PAGES = {xiv+160},
      ISBN = {978-0-691-15371-1},
   MRCLASS = {11G35 (11G15 11G18 11G50 14G05)},
  MRNUMBER = {2918151},
MRREVIEWER = {Yuri\ Bilu},
}

@book {Lang91,
    AUTHOR = {Lang, Serge},
     TITLE = {Number theory. {III}},
    SERIES = {Encyclopaedia of Mathematical Sciences},
    VOLUME = {60},
      NOTE = {Diophantine geometry},
 PUBLISHER = {Springer-Verlag, Berlin},
      YEAR = {1991},
     PAGES = {xiv+296},
      ISBN = {3-540-53004-5},
   MRCLASS = {11Gxx (14Gxx)},
  MRNUMBER = {1112552},
MRREVIEWER = {Joseph\ H.\ Silverman},
       DOI = {10.1007/978-3-642-58227-1},
       URL = {https://doi-org.proxy2.cl.msu.edu/10.1007/978-3-642-58227-1},
}

@article {BL,
    AUTHOR = {Birkenhake, Christina and Lange, Herbert},
     TITLE = {The exponent of an abelian subvariety},
   JOURNAL = {Math. Ann.},
  FJOURNAL = {Mathematische Annalen},
    VOLUME = {290},
      YEAR = {1991},
    NUMBER = {4},
     PAGES = {801--814},
      ISSN = {0025-5831,1432-1807},
   MRCLASS = {14K02 (14K05)},
  MRNUMBER = {1119951},
MRREVIEWER = {Jos\'e\ M.\ Mu\~noz Porras},
       DOI = {10.1007/BF01459272},
       URL = {https://doi-org.proxy2.cl.msu.edu/10.1007/BF01459272},
}

@article {Remond00,
    AUTHOR = {R{\'e}mond, Ga{\"e}l},
     TITLE = {D\'ecompte dans une conjecture de {L}ang},
   JOURNAL = {Invent. Math.},
  FJOURNAL = {Inventiones Mathematicae},
    VOLUME = {142},
      YEAR = {2000},
    NUMBER = {3},
     PAGES = {513--545},
      ISSN = {0020-9910,1432-1297},
   MRCLASS = {14K15 (11G10 14G05)},
  MRNUMBER = {1804159},
MRREVIEWER = {Klaus\ K\"unnemann},
       DOI = {10.1007/s002220000095},
       URL = {https://doi-org.proxy2.cl.msu.edu/10.1007/s002220000095},
}

@article {BeckSchicho,
    AUTHOR = {Beck, Tobias and Schicho, Josef},
     TITLE = {Adjoint computation for hypersurfaces using formal
              desingularizations},
   JOURNAL = {J. Algebra},
  FJOURNAL = {Journal of Algebra},
    VOLUME = {320},
      YEAR = {2008},
    NUMBER = {11},
     PAGES = {3984--3996},
      ISSN = {0021-8693,1090-266X},
   MRCLASS = {14E15 (14F10)},
  MRNUMBER = {2464803},
MRREVIEWER = {D.\ A.\ Stepanov},
       DOI = {10.1016/j.jalgebra.2008.08.002},
       URL = {https://doi-org.proxy2.cl.msu.edu/10.1016/j.jalgebra.2008.08.002},
}

@article {GuerraPirola,
    AUTHOR = {Guerra, Lucio and Pirola, Gian Pietro},
     TITLE = {On the finiteness theorem for rational maps on a variety of
              general type},
   JOURNAL = {Collect. Math.},
  FJOURNAL = {Universitat de Barcelona. Collectanea Mathematica},
    VOLUME = {60},
      YEAR = {2009},
    NUMBER = {3},
     PAGES = {261--276},
      ISSN = {0010-0757,2038-4815},
   MRCLASS = {14E05},
  MRNUMBER = {2560970},
       DOI = {10.1007/BF03191371},
       URL = {https://doi-org.proxy2.cl.msu.edu/10.1007/BF03191371},
}

@article {Jiang,
    AUTHOR = {Jiang, Zhi},
     TITLE = {On {S}everi type inequalities},
   JOURNAL = {Math. Ann.},
  FJOURNAL = {Mathematische Annalen},
    VOLUME = {379},
      YEAR = {2021},
    NUMBER = {1-2},
     PAGES = {133--158},
      ISSN = {0025-5831,1432-1807},
   MRCLASS = {14C20 (14E05)},
  MRNUMBER = {4211084},
MRREVIEWER = {Piotr\ Pokora},
       DOI = {10.1007/s00208-020-02082-6},
       URL = {https://doi-org.proxy2.cl.msu.edu/10.1007/s00208-020-02082-6},
}
\bibliographystyle{amsalpha}

\Addresses

\end{document}